\documentclass{amsart}
\usepackage{amsfonts}

\usepackage{graphicx}

\input{tcilatex}
\input tcilatex

\begin{document}
\title[Population growth under special catastrophic events]{Scaling features of two special Markov chains involving total disasters}
\author{Branda Goncalves, Thierry Huillet}
\address{Laboratoire de Physique Th\'{e}orique et Mod\'{e}lisation\\
CNRS-UMR 8089 et Universit\'{e} de Cergy-Pontoise\\
2 Avenue Adolphe Chauvin, F-95302, Cergy-Pontoise, France\\
Thierry.Huillet@u-cergy.fr}
\maketitle

\begin{abstract}
Catastrophe Markov chain population models have received a lot of attention
in the recent past. We herewith consider two special cases of such models
involving total disasters, both in discrete and in continuous-time.
Depending on the parameters range, the two models can show up a
recurrence/transience transition and, in the critical case, a positive/null
recurrence transition. The collapse transition probabilities are chosen in
such a way that the models are exactly solvable and, in case of positive
recurrence, intimately related to the extended Sibuya and Pareto-Zipf
distributions whose divisibility and self-decomposability properties are
shown relevant. The study includes: existence and shape of the invariant
measure, time-reversal, return time to the origin, contact probability at
the origin, extinction probability, height and length of the excursions, a
renewal approach to the fraction of time spent in the catastrophic state,
scale function, first time to collapse and first-passage times, divisibility
properties.\newline



\end{abstract}

\section{Introduction}

We consider two particular instances of both discrete and continuous-time
Markov chains on the integers subject to state-dependent total disasters
probabilities, as particular cases of similar models with partial
catastrophes. The chosen disaster transition probabilities ensure that the
models are exactly solvable. The Sibuya and Pareto-Zipf distributions and
their relevant divisibility properties are respectively involved in the
analysis. In the discrete-time version of these models, there is a
possibility to move up by one unit with some state-dependent probability and
a complementary collapse probability to return back instantaneously to state
zero (total disaster). Partial reflection at the origin is assumed. The
collapse probability will be a decreasing function of the state, in contrast
with the class of ``house-of-cards'' processes where adding a card to an
already large house is more likely to lead to a breakdown. In both models,
we are able to give necessary and sufficient conditions for the existence
and integrability of the invariant measures and characterize their shape
when they exist. We obtain the precise expression of the laws of the first
time to collapse, first-passage times, return time to the origin and height
and length of the excursions. A renewal approach to the fraction of time
spent in the catastrophic state (zeroset) is also supplied.

The discrete-time models are used as building blocks of their
continuous-time versions, which are obtained after a state-dependent
Poissonization. The jump rates are chosen to be algebraic in the state,
leading again to explicit transience/recurrence and positive/null recurrence
criteria.

\section{A first special Markov chain model with total disaster: results and
background}

Discrete-time integral-valued growth-collapse processes where long periods
of linear growth alternate with rare catastrophic events occur in a large
variety of systems. A collapse or catastrophic event is when the size of the
system shrinks by a random number of units, not exceeding the current
system's size. A total disaster is when size of the system shrinks
instantaneously to zero (a massive extinction event). Disastrous
growth-collapse models occur as models for population growth subject to rare
catastrophic extinction events.\emph{\ }

A one-parameter version of such discrete-time models was investigated in 
\cite{Hui}.\emph{\ }Here, holding probabilities were allowed (with some
probability the system's size can be left unchanged) and pure reflection at
the origin was assumed (once in state zero, the system's size grows by one
unit with probability $1$). Whenever zero is a reflection/absorption
barrier, pomp periods will alternate with periods of scarcity (the Joseph
and Noah effects). We herewith focus on discrete-time disastrous
growth-collapse models with no holding probability and with zero standing
for a reflection/absorption barrier. The probabilities of either growth or
disastrous events will be chosen to be suitably dependent one the current
state, so as to favor large populations in the long run. Also, we shall
consider a continuous time version of this process with jump rates algebraic
in the current state. In this setup, transient and equilibrium issues will
be studied.

We herewith summarize the obtained results for a first Markov chain model
with total disaster, both in discrete and continuous times.

\subsection{\textbf{The model and outline of the results.}}

Our first model can precisely be described as follows.

Let $\beta >0$, $\nu >-1$ and $0<\alpha <\nu +1.$ Consider the discrete
time-homogeneous Markov chain $X:=\left( X_{n};n\geq 0\right) $ with
state-space $\Bbb{N}_{0}=\left\{ 0,1,...\right\} $ and non-homogeneous
spatial transition probabilities characterized by:

$\bullet $ given $X_{n}=x\in \left\{ 1,2,...\right\} $, the increment of $%
X_{n}$ is 
\begin{equation}
\begin{array}{l}
+1\text{ with probability}:\text{ }p_{x}=1-\alpha /\left( \nu +x^{\beta
}\right) \\ 
-x\text{ with probability}:\text{ }q_{x}=\alpha /\left( \nu +x^{\beta
}\right) .
\end{array}
\label{0a}
\end{equation}

$\bullet $ given $X_{n}=0$, the increment of $X_{n}$ is $+1$ with
probability $p_{0}\leq 1$ and $0$ with probability $q_{0}=1-p_{0}.$

This defines the transition matrix $P=\left[ P\left( x,y\right) \right] $ of
the discrete-time (DT) Markov chain $X_{n}$ as: 
\begin{equation*}
P=\left[ 
\begin{array}{llllll}
q_{0} & p_{0} &  &  &  & \cdots \\ 
q_{1} & 0 & p_{1} &  &  & \cdots \\ 
\vdots & \vdots & \ddots & \ddots &  & \cdots \\ 
q_{x} & 0 & \cdots & \ddots & p_{x} & \cdots \\ 
\vdots & 0 & \cdots &  & 0 & \ddots \\ 
\vdots & \vdots &  &  &  & \ddots
\end{array}
\right] .
\end{equation*}
Note that, defining $P^{c}\left( x,y\right) =\sum_{z=0}^{y}P\left(
x,z\right) $, $x^{\prime }>x\Rightarrow P^{c}\left( x^{\prime },y\right)
<P^{c}\left( x,y\right) $ for all $y$ if and only if $q_{x^{\prime
}}<q_{x}<1 $ which is the case when $q_{x}=\alpha /\left( \nu +x^{\beta
}\right) $. The chain $X$ is stochastically monotone. \newline

With $\left( U_{n},\text{ }n\geq 1\right) $ a sequence of independent
identically distributed (iid) uniform random variables (rvs), the dynamics
of $X_{n}$ reads 
\begin{equation*}
X_{n+1}=\left( X_{n}+1\right) \mathbf{1}\left( U_{n+1}>q_{X_{n}}\right) .
\end{equation*}
With $x\geq 1$, we have: 
\begin{eqnarray*}
f\left( x\right)  &=&\Bbb{E}\left( X_{n+1}-X_{n}\mid X_{n}=x\right)
=1-\alpha /\left( \nu +x^{\beta }\right) -\left( \alpha x\right) /\left( \nu
+x^{\beta }\right)  \\
\sigma ^{2}\left( x\right)  &=&\Bbb{E}\left( \left( X_{n+1}-X_{n}\right)
^{2}\mid X_{n}=x\right) =1-\alpha /\left( \nu +x^{\beta }\right) +\left(
\alpha x^{2}\right) /\left( \nu +x^{\beta }\right)  \\
\frac{f\left( x\right) }{\sigma ^{2}\left( x\right) } &=&\frac{\nu -\alpha
+x^{\beta }-\alpha x}{\nu -\alpha +x^{\beta }+\alpha x^{2}}<1
\end{eqnarray*}
giving the local drift and variance of $X_{n}$. Note that, when $x$ is
large, $f\left( x\right) $ is convex (concave) when $\beta <1$ ($\beta \geq 1
$), with: 
\begin{eqnarray*}
f\left( x\right)  &\sim &-\alpha x^{1-\beta }\rightarrow -\infty \text{ when 
}\beta <1 \\
f\left( x\right)  &\sim &1-\alpha x^{-\left( \beta -1\right) }\rightarrow 1%
\text{ when }\beta >1 \\
f\left( x\right)  &=&1-\alpha \frac{1+x}{\nu +x}\rightarrow 1-\alpha \text{
when }\beta =1.
\end{eqnarray*}
When $x$ is large also, 
\begin{eqnarray*}
\frac{f\left( x\right) }{\sigma ^{2}\left( x\right) } &\sim &-\frac{1}{x}%
\text{ if }\beta <1\text{; }\sim \frac{1-\alpha }{1+\alpha }\frac{1}{x}\text{
if }\beta =1 \\
&\sim &\frac{1}{x^{2-\beta }}\text{ if }2>\beta >1\text{; }\sim \frac{1}{%
1+\alpha }\text{ if }\beta =2 \\
&\sim &1^{-}\text{ if }\beta >2.
\end{eqnarray*}
In model (\ref{0a}), the walker $X_{n}$ is occasionally bounced back to the
origin and the probability of this event becomes very small once the walker
has already reached a large value $x$. We will check that:

- If $\beta >1$, the chain is transient with no non-trivial ($\neq \mathbf{0}
$) invariant measure.

- If $0<\beta <1$, the chain is positive recurrent, the invariant measure of
which has stretched exponential behaviour.

- If $\beta =1$ (critical case). When $x$ is large, the drift of this MC is
of order $1-\alpha +\left( \nu -\alpha \right) /x$. So when $\alpha >1$, the
walker is attracted to the origin: The strength of the attraction goes like $%
\alpha -1$ for large $x$. For $\alpha <1$, the walker is repelled from the
origin correspondingly. When $\alpha =1$, its drift is still attracting but
of order $\left( \nu -1\right) /x$, and the drift that the walker feels
vanishes when $x$ approaches $\infty $ if $\nu <1$. We will see that while
crossing the critical value $\alpha =1$ from above, the process $X$ switches
from positive recurrent to null-recurrent.\newline

\textbf{Remark:} The collapse probability $q_{x}$ may be a decreasing or an
increasing function of the current population size $x.$ In the former case,
large populations are getting more and more stable, after having survived
the early growth stages. This is the case under study here in model (\ref{0a}%
) with $q_{x}\underset{x\rightarrow \infty }{\sim }\alpha x^{-\beta }$, ($%
\alpha ,\beta >0$). In the latter opposite case, large populations would be
more susceptible and vulnerable to collapse and so quite unlikely to grow
large and develop too much. This would be the case for a model with $q_{x}%
\underset{x\rightarrow \infty }{\sim }1-\alpha x^{-\beta }$ (while switching
the role of $p_{x}$ and $q_{x}$). If one also thinks of the process of
building a house of cards, clearly adding a new card to a house of cards of
size $x$ is more likely to lead to a collapse of the whole structure if $x$
is already large (\cite{FFG}, p. $47$ and \cite{FG}, p. $9$). In such a
situation, the process is always positive recurrent with light-tailed
invariant measure (stretched exponential). $\Box $\newline

We now come to a natural continuous-time version of the process (\ref{0a})
which will be one of our main process of interest.

Let $\lambda \in \left( -\infty ,+\infty \right) $ and consider the
transition rate matrix $Q$ of a continuous-time MC process $\overline{X}%
\left( t\right) :$%
\begin{equation}
Q=D_{\mathbf{r}}\left( P-I\right)  \label{0b}
\end{equation}
with $D_{\mathbf{r}}=$diag$\left( \mathbf{r}\right) $ the diagonal matrix
formed from the rate vector $\mathbf{r}=(r_{0},r_{1},...,r_{x},...)$ with $%
r_{x}=r_{0}\left( x+1\right) ^{\lambda }$, $x\geq 0$, $r_{0}>0$. Clearly, $%
X_{n}$ is the embedded MC of $\overline{X}\left( t\right) $.

For such a continuous-time model, the growth rate of the transition $%
x\rightarrow x+1$ is $r_{x}p_{x}$ with $r_{x}p_{x}\sim r_{0}x^{\lambda }$
for large $x$, while the one of the collapse transition $x\rightarrow 0$ is $%
r_{x}q_{x}$ with $r_{x}q_{x}\sim \alpha r_{0}x^{\lambda -\beta }$ for large $%
x.$ If $\lambda >0,$ $\alpha x^{\lambda -\beta }\ll x^{\lambda }$ always and
if $0<\lambda <\beta $, $\overline{X}$ moves up by one unit frequently,
while its collapse becomes increasingly rare. If $\lambda <0$ (the lazy
chain), the growth rate $r_{x}p_{x}\sim r_{0}x^{\lambda }$ is small for
large $x$ while the collapse rate $r_{x}q_{x}\sim \alpha r_{0}x^{\lambda
-\beta }$ is still smaller. In any case, the collapse rates are small
compared to the growth ones.

We can view the process $\overline{X}\left( t\right) $ as follows. Let $%
P\left( t\right) $ be a standard Poisson process with intensity $t\geq 0$
and $P\left( 0\right) =0.$ Let $Z\left( t\right) =X_{P\left( t\right) }$ be
the chain $X_{n}$ subordinated to $P\left( t\right) $. Clearly $Z\left(
t\right) $ is a CT Markov chain with transition rate matrix $P-I.$ Then the
process $\overline{X}\left( t\right) $ turns out to be 
\begin{equation*}
\overline{X}\left( t\right) =Z\left( \int_{0}^{t}r_{\overline{X}%
_{s}}ds\right) .
\end{equation*}
It has infinitesimal backward generator $G_{\overline{X}}$ whose action on
real-valued bounded functions $h$ on $\Bbb{N}_{0}=\left\{ 0,1,2,...\right\} $
is 
\begin{equation*}
G_{\overline{X}}h\left( x\right) =r_{x}\left\{ \left( \left( h\left(
x+1\right) -h\left( x\right) \right) p_{x}+\left( h\left( 0\right) -h\left(
x\right) \right) q_{x}\right) \right\} ,
\end{equation*}
meaning 
\begin{equation*}
\Bbb{E}_{x}h\left( \overline{X}\left( t\right) \right) =\Bbb{E}_{x}h\left( 
\overline{X}\left( 0\right) \right) +\int_{0}^{t}\Bbb{E}_{x}\left( G_{%
\overline{X}}h\right) \left( \overline{X}\left( s\right) \right) ds.
\end{equation*}
When $h\left( x\right) =x$, with $x\left( t\right) :=\Bbb{E}_{x}\overline{X}%
\left( t\right) $, $x\left( 0\right) =x$, by Jensen inequality 
\begin{eqnarray*}
\overset{.}{x}\left( t\right) &=&\Bbb{E}_{x}\left( r\left( \overline{X}%
\left( t\right) \right) \left( p_{\overline{X}\left( t\right) }-\overline{X}%
\left( t\right) q_{\overline{X}\left( t\right) }\right) \right) =:\Bbb{E}_{x}%
\overline{f}\left( \overline{X}\left( t\right) \right) \\
&=&r_{0}\Bbb{E}_{x}\left( \left( \overline{X}\left( t\right) +1\right)
^{\lambda }\left( \frac{\nu -\alpha \left( 1+\overline{X}\left( t\right)
\right) +\overline{X}\left( t\right) ^{\beta }}{\nu +\overline{X}\left(
t\right) ^{\beta }}\right) \right) \\
&\geq &\text{ }r_{0}\left( x+1\right) ^{\lambda }\left( \frac{\nu -\alpha
\left( 1+x\right) +x^{\beta }}{\nu +x^{\beta }}\right) =r_{0}\left(
x+1\right) ^{\lambda }f\left( x\right)
\end{eqnarray*}
if the latter function to the right of the inequality (the drift $\overline{f%
}\left( x\right) =r_{0}\left( x+1\right) ^{\lambda }f\left( x\right) $ of
the continuous-time MC) is a convex function of $x$, at least for large $x$.

The range $\lambda >0$ ($\lambda <0$) accounts for the fact that the moves
of $\overline{X}$ get frequent (respectively rare) when the height $x$ of $%
\overline{X}$ gets large, and given a move has occurred and $x$ is large, $%
\overline{X}$ grows by one unit with large probability $\sim 1-\alpha
x^{-\beta }$ or undergoes a catastrophic event with small complementary
probability $\sim \alpha x^{-\beta }$. Such transition mechanisms favor
large values of $\overline{X}$. This chain is irreducible and aperiodic,
either transient or recurrent (possibly then either positive or null
recurrent). We will show that:

$\left( i\right) $ When it is transient ($\beta >1$), the process $\overline{%
X}$ is either explosive ($\lambda >1$) or non-explosive ($\lambda \leq 1$).
When $\beta >1$, after a finite number of returns to $0$, $\overline{X}$
drifts to $\infty $. And (from an argument on Yule processes to appear in
the proof) it explodes if and only if $\lambda >1.$ The process $\overline{X}
$ has no non-trivial ($\neq \mathbf{0}$) invariant measure.

$\left( ii\right) $ When it is recurrent ($\beta \leq 1$), the process $%
\overline{X}$ is:

- recurrent positive if $\beta <1$, with invariant measure showing a
stretched exponential behavior.

- recurrent positive if $\beta =1$ and $\alpha +\lambda >1$, recurrent null
if $\beta =1$ and $\alpha +\lambda \leq 1.$ The invariant measure of $%
\overline{X}$ is of power-law type with index $\alpha +\lambda $. The
critical value $\beta =1$ separates a recurrent phase ($\beta <1$) from a
transient phase ($\beta >1$).

The birth and collapse probabilities depend on the current state of the
population in the specific way just described, together with the jump rates.
Our study will include, among other topics, first-return time probabilities
to the origin (excursion length), eventual return (contact) probability to
the origin, excursion height, time to failure and the fraction of time spent
in the catastrophic state.

The continuous-time version of model (\ref{0a}) is skip-free to the right,
with $p_{x},$ $q_{x}$ and $r_{x}$ dependent on $x$ as specified. To a large
extent therefore, the Markov chain model under study here is one of the
simplest possible in the vast family of growth-collapse models. As we will
show, it turns out that it is an exactly solvable case. Using ideas stemming
from excursion theory, we will make precise here to what extent growth is
(algebraically) slow when catastrophic events are rare under our
over-simplistic model hypothesis. It is hoped that models in the same class
of universality could share similar behaviors as one reasonably can expect
some sort of robustness. This will be the purpose of the last Section $4$
involving a Pareto-Zipf paradigm.

\section{Detailed analysis of the special Markov chain with total disaster:
proofs}

In this Section $3$, we will supply a detailed study of the general DT
catastrophe Markov chain (MCC), without specifying the disaster
probabilities $q_{x}.$ For each item under study, we will subsequently
particularize the detailed expression of the results for the special MC (\ref
{0a}). We end up this Section with the study of its continuous-time version
with algebraic rates $r_{x}.$

\subsection{\textbf{The discrete-time MCC chain}}

Consider a general catastrophe process $X_{n}$ for which both $p_{x}$ and $%
q_{x}>0$, for all $x\geq 0$, with $p_{x}+$ $q_{x}=1$ and so with associated
stochastic transition matrix: $P=\left[ P\left( x,y\right) \right] $, $%
\left( x,y\right) \in \Bbb{N}_{0}^{2}$ with $P\left( x,0\right) =q_{x}$ and $%
P\left( x,x+1\right) =p_{x},$ $x\geq 0.$

\subsubsection{\textbf{Existence and shape of the invariant measure. }}

Let $\mathbf{\pi }^{\prime }\equiv \left( \pi _{0},\pi _{1},..\right) $ be
the row-vector of the invariant measure, whenever it exists. Then $\mathbf{%
\pi }$ should solve $\mathbf{\pi }^{\prime }=\mathbf{\pi }^{\prime }P,$
whose formal solution is: 
\begin{equation}
p_{0}\pi _{0}=\sum_{x\geq 1}\pi _{x}q_{x}\text{ and }\pi _{x}=\pi
_{0}\prod_{y=0}^{x-1}p_{y}\text{, }x\geq 1.  \label{1}
\end{equation}
Let $u_{x}=\prod_{y=1}^{x-1}p_{y}$ with $u_{1}\equiv 1$. Using the second
equation, the first equation is satisfied whenever 
\begin{equation*}
\sum_{x\geq 1}q_{x}\prod_{y=1}^{x-1}p_{y}=\sum_{x\geq 1}\left(
u_{x}-u_{x+1}\right) =1,
\end{equation*}
so also when $u_{\infty }=\prod_{y=1}^{\infty }p_{y}=0$ which is fulfilled
if and only if 
\begin{equation*}
C_{1}\equiv \sum_{y=1}^{\infty }q_{y}=\infty .
\end{equation*}

We first conclude that there exists an invariant measure if and only if $%
C_{1}=\infty .$

If in addition, 
\begin{equation*}
C_{2}\equiv \sum_{x\geq 1}\prod_{y=0}^{x-1}p_{y}<\infty ,
\end{equation*}
then $\pi _{0}=\frac{1}{1+C_{2}}\in \left( 0,1\right) $ and the invariant
measure is unique and is a proper invariant probability measure$.$ In this
case, with the empty product being $1$, we have 
\begin{equation*}
\pi _{x}=\frac{\prod_{y=0}^{x-1}p_{y}}{1+C_{2}}\text{, }x\geq 0.
\end{equation*}

When $C_{2}=\infty $, the measure $\pi _{x}$ solution to (\ref{1}) exists
but it is not a probability measure as there is no way to normalize it so as
to have a probability measure. We will see below that this corresponds to a
case where the first return time to $0$\ of the chain has infinite mean.%
\newline

APPLICATION: When dealing with the special MC, we conclude:

- $C_{1}<\infty $ if and only if $\beta >1$: in this case the MC is
transient with no invariant measure.

- If $\beta <1,$ then $C_{1}=\infty $ and $C_{2}<\infty :$ the MC is
positive recurrent. Furthermore, for large $x$: 
\begin{equation*}
\pi _{x}\sim \prod_{y=0}^{x-1}\left( 1-\frac{\alpha }{\nu +y^{\beta }}%
\right) \sim e^{-\alpha \sum^{x}y^{-\beta }}\sim e^{-\frac{\alpha }{1-\beta }%
x^{1-\beta }}
\end{equation*}
with stretched exponential behaviour.\newline

\textbf{Remark:}\emph{\ }the simpler case $\beta =0$\ was excluded from the
study. Here, the transition probabilities $\left( p_{x},q_{x}\right) $, $%
x\geq 1,$\ are homogeneous and may be set to $\left( p,q\right) $, $p+q=1,$
where $q=\alpha /\left( \nu +1\right) <1.$\ The invariant measure in this
case exists, is summable and takes the simple geometric form $\pi _{x}=\pi
_{0}p_{0}p^{x-1}$\ if $x\geq 1.$\ The value of $\pi _{0}=1/\left(
1+C_{2}\right) $\ is found to be: $\pi _{0}=q/\left( p_{0}+q\right) .$\ The
corresponding chain is positive recurrent. $\Box $

- In the critical case $\beta =1$: 
\begin{equation*}
\pi _{x}\sim \prod_{y=0}^{x-1}\left( 1-\frac{\alpha }{\nu +y}\right) \sim
e^{-\alpha \sum^{x}y^{-1}}\sim x^{-\alpha }
\end{equation*}
with power-law$\left( \alpha \right) $ behaviour. The chain is recurrent. We
will compute the exact shape of $\pi _{x}$ later. But from this one can
guess that the DT chain is positive recurrent if $\alpha >1$, null recurrent
if $\alpha \leq 1.$

\subsubsection{\textbf{Return time to the origin.}}

Let $X_{n}$ be the Markov chain with transition probability matrix $P$ and
state-space $\Bbb{N}_{0}.$ Starting from $X_{0}=x\geq 1$, the walker moves
one step up with probability $p_{x}$ or goes back to the origin (the
catastrophic event) with probability $q_{x}$. Once at $0$, the walker moves
one step up with probability $p_{0}$ or stays at $0$ with probability $q_{0}$
(the latter event being considered as a one-step return to $0$). Clearly of
interest are the times $\tau _{0,0}\geq 1$ between consecutive visits to $0$
(the first return times to $0$). With the convention that the empty product
is $1$, with $u_{0}\equiv 1$, we have 
\begin{equation}
\Bbb{P}\left( \tau _{0,0}=x+1\right)
=q_{x}\prod_{y=0}^{x-1}p_{y}=p_{0}\left( u_{x}-u_{x+1}\right) \text{, }x\geq
0.  \label{1a}
\end{equation}
Equivalently, $\Bbb{P}\left( \tau _{0,0}>x\right)
=\prod_{y=0}^{x-1}p_{y}=p_{0}u_{x}.$ We note that, conventionally, $\tau
_{0,0}=1$ with probability $q_{0}$ (the holding probability at state $0$)$.$

Note also that $p_{x}=\Bbb{P}\left( \tau _{0,0}>x+1\right) /\Bbb{P}\left(
\tau _{0,0}>x\right) $ and $q_{x}=\Bbb{P}\left( \tau _{0,0}=x+1\right) /\Bbb{%
P}\left( \tau _{0,0}>x\right) :$ If the law of the lifetime $\tau _{0,0}$ is
known in the first place, this gives the survival probability $p_{x}$ given
the age (backward recurrence time) $X_{n}$ of the current machine is $x.$

Note also that $\Bbb{P}\left( \tau _{0,0}<\infty \right) =1$ if and only if $%
u_{\infty }=0$ ($C_{1}=\infty $) which is the recurrence condition for $%
X_{n} $.

We conclude that: If $C_{1}\equiv \sum_{y=1}^{\infty }q_{y}<\infty ,$ $X_{n}$
is transient with $\Bbb{P}\left( \tau _{0,0}<\infty \right) <1$. If $%
C_{1}=\infty ,$ $X_{n}$ is recurrent. If $C_{1}=\infty $ and $C_{2}\equiv
\sum_{x\geq 1}\prod_{y=0}^{x-1}p_{y}<\infty ,$ $X_{n}$ is positive recurrent
with $\mu :=\Bbb{E}\left( \tau _{0,0}\right) =1/\pi _{0}=1+C_{2}<\infty .$
If $C_{1}=C_{2}=\infty ,$ $X_{n}$ is null recurrent with $\tau _{0,0}<\infty 
$ almost surely (a.s.) and $\Bbb{E}\left( \tau _{0,0}\right) =\infty $. Note
that in the positive recurrent case, 
\begin{equation*}
\pi _{x}=\frac{\prod_{y=0}^{x-1}p_{y}}{1+C_{2}}=\frac{\Bbb{P}\left( \tau
_{0,0}>x\right) }{\mu }\text{, }x\geq 0.
\end{equation*}
\newline
APPLICATION:

- If $\beta >1:$ the chain is transient with ($0<u_{\infty }<1$) 
\begin{equation*}
\Bbb{P}\left( \tau _{0,0}=\infty \right) =u_{\infty }=\prod_{y=0}^{\infty
}\left( 1-\frac{\alpha }{\nu +y^{\beta }}\right) >0.
\end{equation*}

- If $\beta <1:$ the chain is positive recurrent with 
\begin{equation*}
\mu :=\Bbb{E}\left( \tau _{0,0}\right) =1+C_{2}=1+\sum_{x\geq
1}\prod_{y=0}^{x-1}\left( 1-\frac{\alpha }{\nu +y^{\beta }}\right) <\infty .
\end{equation*}

- Critical case: $\beta =1$, $0<\alpha <\nu +1$. In this case, let 
\begin{equation*}
\psi _{0}\left( z\right) :=\sum_{x\geq
1}q_{x}z^{x}\prod_{y=1}^{x-1}p_{y}=\sum_{x\geq 1}z^{x}\frac{\alpha }{\nu +x}%
\prod_{y=1}^{x-1}\left( 1-\frac{\alpha }{\nu +y}\right) =\frac{\alpha z}{\nu
+1}\cdot F\left( 1,\nu +1-\alpha ;\nu +2;z\right)
\end{equation*}
involving a Gauss hypergeometric function (by $F$ we mean here $_{2}F_{1}$).
The function $\psi _{0}\left( z\right) $ is the probability generating
function (pgf) of an extended Sibuya rv $S\geq 1$ with parameters $\left(
\alpha ,\nu \right) ,$ \cite{KP}. The shifted rv $S-1$ (with pgf $z^{-1}\psi
_{0}\left( z\right) $) is discrete-self-decomposable (SD), \cite{KP} (see
the Appendix $8.1$ for a reminder on this concept). When $\nu =0$, the law
of $S$ is also known as the standard Sibuya distribution, \cite{Sibu}, while
when $\nu =\alpha $ the Yule-Simon distribution is recovered, \cite{Yule}, 
\cite{Sim}.

When $\alpha =1$, an easy computation shows that 
\begin{equation*}
\psi _{0}\left( z\right) =z-\nu L_{\nu }\left( z\right) \left( 1-z\right)
\end{equation*}
where 
\begin{equation*}
L_{\nu }\left( z\right) =\sum_{n\geq 1}z^{n}/\left( \nu +n\right) \underset{%
z\downarrow 1}{\sim }-\log \left( 1-z\right) .
\end{equation*}

Using Stirling formula: 
\begin{eqnarray*}
\Bbb{P}\left( S=x\right) &:&=s_{x}=\left[ z^{x}\right] \psi _{0}\left(
z\right) \underset{x\rightarrow \infty }{\sim }\alpha \Gamma \left( \nu
+1\right) x^{-\left( \alpha +1\right) }/\Gamma \left( \nu +1-\alpha \right) 
\text{ if }\alpha \neq 1 \\
\Bbb{P}\left( S>x\right) &=&\frac{1}{\nu +1}\frac{1}{\nu +x}\text{ if }%
\alpha =1\text{ (}\nu >0\text{)}
\end{eqnarray*}
and for all $x\geq 1$%
\begin{equation}
s_{x+1}/s_{x}=\left( \nu -\alpha +x\right) /\left( \nu +x+1\right) <1,\text{ 
}s_{1}=\alpha /\left( \nu +1\right)  \label{f13a}
\end{equation}
($s_{x}$ is monotone decreasing). As a result, $s_{1}$ is the maximal value
of the $s_{x}$: the probability mass function (pmf) $s_{x}$ has its mode at $%
x=1.$

Equivalently, 
\begin{equation*}
\begin{array}{l}
\psi _{0}\left( z\right) \underset{z\downarrow 1}{\sim }1-\frac{\Gamma
\left( \nu +1\right) \Gamma \left( 1+\alpha \right) }{\Gamma \left( \nu
+1-\alpha \right) }\left( 1-z\right) ^{\alpha }\text{ if }0<\alpha <1\text{ }
\\ 
\psi _{0}\left( z\right) \underset{z\downarrow 1}{\sim }1-\frac{\nu }{\alpha
-1}\left( 1-z\right) +O\left( 1-z\right) ^{\alpha }\text{ if }\alpha >1 \\ 
\psi _{0}\left( z\right) \underset{z\downarrow 1}{\sim }1-\left( 1-z\right)
-\nu \left( 1-z\right) \log \left( 1-z\right) \text{ if }\alpha =1
\end{array}
\end{equation*}
The pgf of $\tau _{0,0}$ itself therefore reads,

\begin{eqnarray*}
\phi _{0,0}\left( z\right) &=&\sum_{x\geq
0}q_{x}z^{x+1}\prod_{y=0}^{x-1}p_{y}=z\left( q_{0}+p_{0}\sum_{x\geq
1}q_{x}z^{x}\prod_{y=1}^{x-1}p_{y}\right) \\
&=&z\left( q_{0}+p_{0}\frac{\alpha z}{\nu +1}\cdot F\left( 1,\nu +1-\alpha
;\nu +2;z\right) \right) =z\left( q_{0}+p_{0}\psi _{0}\left( z\right)
\right) .
\end{eqnarray*}
Using $F\left( a,b;c;1\right) =\frac{\Gamma \left( c\right) \Gamma \left(
c-a-b\right) }{\Gamma \left( c-a\right) \Gamma \left( c-b\right) }$\ and $%
F^{\prime }\left( a,b;c;1\right) =\frac{ab}{c}F\left( a+1,b+1;c+1;1\right) $%
\ we find the mean persistence time 
\begin{equation*}
\phi _{0,0}^{\prime }\left( 1\right) =\Bbb{E}\left( \tau _{0,0}\right) :=\mu
=\frac{1}{\pi _{0}}=1+p_{0}\frac{\nu }{\alpha -1}.
\end{equation*}
It can be checked that if $\Bbb{E}\left( \tau _{0,0}\right) $\ is to exist,
then necessarily $\Bbb{E}\left( \tau _{0,0}\right) >1+p_{0}$.

* This condition forces $1<\alpha <\nu +1$\ which is the positive recurrence
condition for the critical MCC. If $p_{0}=1$, $\Bbb{E}\left( \tau
_{0,0}\right) =\frac{\nu +\alpha -1}{\alpha -1}>2$ and there are no trivial
excursions.\newline

In this positive recurrent case, we also have

\begin{eqnarray*}
\pi _{x} &=&\frac{\prod_{y=0}^{x-1}p_{y}}{\mu }=\pi
_{0}p_{0}\prod_{y=1}^{x-1}\left( 1-\frac{\alpha }{\nu +y}\right) \text{, }%
x\geq 0 \\
\frac{\pi _{x}}{1-\pi _{0}} &=&\frac{\pi _{0}p_{0}}{1-\pi _{0}}%
\prod_{y=1}^{x-1}\left( 1-\frac{\alpha }{\nu +y}\right) =\frac{\alpha -1}{%
\nu }\prod_{y=1}^{x-1}\left( 1-\frac{\alpha }{\nu +y}\right) \\
&=&\frac{\alpha -1}{\nu -1+x}\prod_{y=1}^{x-1}\left( 1-\frac{\alpha -1}{\nu
-1+y}\right) \text{, }x\geq 1
\end{eqnarray*}
an identity showing that $Y_{\infty }:=X_{\infty }\mid X_{\infty }\geq 1%
\overset{d}{\sim }$ Sibuya$\left( \alpha -1,\nu -1\right) $\ so with 
\begin{eqnarray*}
\Bbb{E}\left( z^{Y_{\infty }}\right) &=&\frac{\left( \alpha -1\right) z}{\nu 
}\cdot F\left( 1,\nu +1-\alpha ;\nu +1;z\right) =:\psi _{\infty }\left(
z\right) \\
\Bbb{E}\left( z^{X_{\infty }}\right) &=&\pi _{0}+\left( 1-\pi _{0}\right)
\psi _{\infty }\left( z\right) .
\end{eqnarray*}
This gives an explicit expression of the pgf of $X_{\infty }$ (in the
positive recurrent case) in terms of a Gauss hypergeometric function. Note
that $X_{\infty }$ is the Bernoulli mixture of the two rvs $Y_{\infty
}:=X_{\infty }\mid X_{\infty }\geq 1$ and (say) $Y_{0}$ which is degenerate
at $0$ so with $Y_{0}\overset{d}{\sim }\delta _{0}.$ We also observe that $%
\pi _{x}$ is unimodal with mode at the origin ($\pi _{x+1}/\pi _{x}=\left( 1-%
\frac{\alpha }{\nu +x}\right) <1$, $x\geq 1$, $\pi _{1}/\pi _{0}=p_{0}<1$).

As an extended Sibuya$\left( \alpha -1,\nu -1\right) $ rv, the shifted
Sibuya rv 
\begin{equation*}
Y_{\infty }-1:=\left( X_{\infty }-1\mid X_{\infty }\geq 1\right) \text{ is
discrete self-decomposable.}
\end{equation*}

We now raise the following question: Is $X_{\infty }$ itself
infinitely-divisible, say ID (meaning compound-Poisson)? discrete
self-decomposable? Discrete-SD rvs constitute a remarkable sub-class of ID
rvs, \cite{SH}; see the Appendix $8.1$. It turns out that $X_{\infty }$ is
ID if $p_{0}$ is small enough (else $\pi _{0}$ large enough).

Observing $\pi _{x}^{2}\leq \pi _{x-1}\pi _{x+1}$ for $x\geq 2$,\emph{\ }a
sufficient condition for\emph{\ }$X_{\infty }$ to be log-convex and so ID
(see \cite{SH} Example $11.9$, page $84$) is that 
\begin{equation*}
\pi _{1}^{2}\leq \pi _{0}\pi _{2}\text{ which is }p_{0}\leq 1-\alpha /\left(
\nu +1\right) .
\end{equation*}
But this condition is not necessary although, from the unimodality of $\pi
_{x}$, it could be a necessary and sufficient condition for $X_{\infty }$ to
be self-decomposable in some cases. Whenever $X_{\infty }$ is ID or SD, it
can be produced as limit laws of alternative Markov processes (see Appendix $%
8.2$). We illustrate this result by the following explicit example showing
that this can happen. We were not able to treat the general case. \newline

\textbf{Example:} Suppose $\nu =1$, so with $\pi _{0}=1-p_{0}/\left( \alpha
-1+p_{0}\right) $. Then, with $\alpha \in \left( 1,2\right) :$ 
\begin{eqnarray*}
\psi _{\infty }\left( z\right) &=&\left( \alpha -1\right) zF\left(
1,2-\alpha ;2;z\right) =1-\left( 1-z\right) ^{\alpha -1} \\
\Bbb{E}\left( z^{X_{\infty }}\right) &=&1-\left( 1-\pi _{0}\right) \left(
1-z\right) ^{\alpha -1}
\end{eqnarray*}
The rv $X_{\infty }$ is obtained while thinning (scaling) a Sibuya$\left(
\alpha -1\right) $ distributed rv (with pgf $\psi _{\infty }\left( z\right) $%
), \cite{Sch}. Note that $X_{\infty }$ has infinite mean. Furthermore, (see
Th. $1$ and $2$ of \cite{Chris1}):

- $X_{\infty }$ is ID if and only if: $1-\pi _{0}\leq 1-\left( \alpha
-1\right) =2-\alpha ,$ else: $p_{0}\leq 2-\alpha $ (compare with the
sufficient condition which here is $p_{0}\leq 1-\alpha /2$)$.$

- $X_{\infty }$ is discrete-SD if and only if: $1-\pi _{0}\leq \frac{%
1-\left( \alpha -1\right) }{1+\left( \alpha -1\right) }=\frac{2-\alpha }{%
\alpha },$ else: $p_{0}\leq 1-\alpha /2$ (the right condition).

Note that in this example 
\begin{eqnarray*}
\psi _{0}\left( z\right) &=&\frac{\alpha z}{2}\cdot F\left( 1,2-\alpha
;3;z\right) =\frac{1}{\alpha -1}\left( \alpha -\frac{1}{z}\left( 1-\left(
1-z\right) ^{\alpha }\right) \right) \\
\phi _{0,0}\left( z\right) &=&z\left( q_{0}+p_{0}\psi _{0}\left( z\right)
\right) =z\left( 1+\frac{p_{0}}{\alpha -1}\right) -\frac{p_{0}}{\alpha -1}%
\left( 1-\left( 1-z\right) ^{\alpha }\right) \\
&=&1-\mu \left( 1-z\right) +\left( \mu -1\right) \left( 1-z\right) ^{\alpha }
\end{eqnarray*}
with $\tau _{0,0}$ having finite mean but infinite variance. $\Box $\newline

* If the condition $1<\alpha <\nu +1$ is not met, the special MCC is null
recurrent with $\Bbb{E}\left( \tau _{0,0}\right) =\infty $.

\subsubsection{\textbf{Relation to renewal theory (in the positive recurrent
case). }}

There is a clear connection of the catastrophe model with the discrete
theory of renewal processes generated by $\tau _{0,0}$ [as from \cite{SH}
pages 462-463].

Let the instant $n$ be such that $n\geq X_{n}>0.$ Define $\tau _{n}^{B}$ the
(backward) time elapsed since the last visit of $X_{n}$ to $0$ and $\tau
_{n}^{F}$ the (forward) time to its next visit to $0$. If $n$ is such that $%
X_{n}=0$, we set conventionally $\tau _{n}^{B}=0.$ Thus for each $n$, $\tau
_{n}^{B}\geq 0$, $\tau _{n}^{F}\geq 1.$ Define $\tau _{n}=\tau _{n}^{B}+\tau
_{n}^{F}\geq 1$ the length of the temporal interval between consecutive
visits to $0$ to which $n$ belongs (the length of its attached excursion).

When $\mu <\infty $, both $\tau _{n}^{B}$ and $\tau _{n}^{F}$ (and so $\tau
_{n}$) converge in distribution as $n\rightarrow \infty $ to some limiting
rvs $\tau _{\infty }^{B}$ and $\tau _{\infty }^{F}$ (respectively $\tau
_{\infty }$) with 
\begin{eqnarray*}
\Bbb{P}\left( \tau _{\infty }^{B}=x\text{, }\tau _{\infty }^{F}-1=x^{\prime
}\right) &=&\frac{1}{\mu }\Bbb{P}\left( \tau _{0,0}=x+x^{\prime }+1\right) 
\text{, }x,x^{\prime }\geq 0, \\
\Bbb{P}\left( \tau _{\infty }=x\right) &=&\frac{x\Bbb{P}\left( \tau
_{0,0}=x\right) }{\mu }\text{, }x\geq 1.
\end{eqnarray*}
The law of $\tau _{\infty }$ is the size-biased version of the law of $\tau
_{0,0}$. So $\tau _{n}^{B}\overset{d}{=}\tau _{n}^{F}-1$ with 
\begin{equation*}
\Bbb{P}\left( \tau _{\infty }^{B}=x\right) =\frac{1}{\mu }\Bbb{P}\left( \tau
_{0,0}>x\right) \text{, }x\geq 0.
\end{equation*}
In fact, it can be checked that (see e.g. \cite{SH}, Lemma $9$, p. $447$) 
\begin{equation*}
\left( \tau _{\infty }^{B}\text{, }\tau _{\infty }^{F}-1\right) \overset{d}{=%
}\left( U\circ \left( \tau _{\infty }-1\right) ,\left( 1-U\right) \circ
\left( \tau _{\infty }-1\right) \right) .
\end{equation*}

We conclude that $\pi _{x}$ is the law of the limiting ($n\rightarrow \infty 
$) backward recurrence time $\tau _{\infty }^{B}$ of the discrete renewal
process generated by $\tau _{0,0}\geq 1$.

In the latter displayed formula, with $U$ a uniform random variable
independent of $\tau $, $U\circ \tau _{\infty }$ is the $U-$thinning of $%
\tau _{\infty }:$ $U\circ \tau _{\infty }=\sum_{l=1}^{\tau _{\infty
}}B_{l}\left( U\right) $ where, given $U,$ $B_{l}\left( U\right) ;l\geq 1$
are mutually independent and independent of $\tau _{\infty }$ rvs with law
Bernoulli$\left( U\right) .$ Then

\begin{equation*}
\Bbb{P}\left( U\circ \left( \tau _{\infty }-1\right) =x\right) =\int_{0}^{1}%
\Bbb{P}\left( u\circ \left( \tau _{\infty }-1\right) =x\right) du
\end{equation*}
and 
\begin{equation*}
\Bbb{P}\left( u\circ \left( \tau _{\infty }-1\right) =x\right) =\sum_{y\geq
x}\binom{y}{x}u^{x}\left( 1-u\right) ^{y-x}\Bbb{P}\left( \tau _{\infty
}=y+1\right) .
\end{equation*}
Thus, 
\begin{eqnarray*}
\Bbb{P}\left( U\circ \left( \tau _{\infty }-1\right) =x\right)
&=&\sum_{y\geq x}\binom{y}{x}\int_{0}^{1}u^{x}\left( 1-u\right) ^{y-x}du\Bbb{%
P}\left( \tau _{\infty }=y+1\right) \\
&=&\sum_{y\geq x}\frac{1}{y+1}\Bbb{P}\left( \tau _{\infty }=y+1\right) =%
\frac{1}{\mu }\Bbb{P}\left( \tau _{0,0}>x\right) =\pi _{x}.
\end{eqnarray*}

\subsubsection{\textbf{The} \textbf{fraction of time spent by }$X$\textbf{\
in the catastrophic state}.}

A `trivial' length $1$ excursion appears whenever for some $n$, $%
X_{n}=X_{n+1}=0$ (an event with probability $q_{0}$). If $X_{n}=0$ and $%
X_{n+1}\neq 0,$ a `true' excursion with length (say $\tau _{0,0}^{+}$) at
least $2$ starts from $n,$ ending up when $X$ first revisits $0$. Consider a
time $n$ for which $X_{n-1}\neq 0$ and $X_{n}=0$. We ask for the
distribution of the time elapsed from $n,$ say $\Delta =N-n$ for which $%
X_{N-1}\neq 0$ and $X_{N}=0$ for the first time $N>n$ again. It holds 
\begin{equation*}
\Delta =\tau _{0,0}^{0}+\tau _{0,0}^{+},
\end{equation*}
where $\tau _{0,0}^{0}\overset{d}{\sim }$Geo$\left( p_{0}\right) $\ (viz $%
\Bbb{P}\left( \tau _{0,0}^{0}=k\right) =p_{0}q_{0}^{k}$, $k\geq 0$) is the
idle period (the number of consecutive trivial excursions) and $\tau
_{0,0}^{+}\geq 2$\ the length of a true excursion with: $\Bbb{P}\left( \tau
_{0,0}^{+}=x+1\right) =q_{x}\prod_{y=1}^{x-1}p_{y}$, $x\geq 1$ (the busy
period). The lengths $\tau _{0,0}^{0}$, $\tau _{0,0}^{+}$ are mutually
independent.

We observe that the random variable $\tau _{0,0}^{+}-1$ (which is also $%
H^{+} $, the height of a true excursion), with support $\Bbb{N,}$ can be
generated as follows: 
\begin{equation}
\tau _{0,0}^{+}-1\overset{d}{=}H^{+}=\inf \left( x\geq 1:B_{x}\left( \alpha
,\nu \right) =1\right) ,  \label{bern2}
\end{equation}
where $\left( B_{x}\left( \alpha ,\nu \right) ;\text{ }x\in \Bbb{N}%
_{0}\right) $ is a sequence of independent Bernoulli rvs obeying $\Bbb{P}%
\left( B_{x}\left( \alpha ,\nu \right) =1\right) =\alpha /\left( \nu
+x\right) =q_{x}.$ It is identified to the rv $S$ with pgf $\psi _{0}\left(
z\right) $ defined above.\newline

APPLICATION (positive recurrent case, $\beta =1$, $\alpha >1$): Clearly,
with $C_{2}^{+}=\sum_{x\geq 1}\prod_{y=1}^{x-1}p_{y}=C_{2}/p_{0}<\infty $%
\begin{equation*}
\Bbb{E}\left( \tau _{0,0}^{0}\right) =p_{0}/q_{0}\text{ and }\Bbb{E}\left(
\tau _{0,0}^{+}\right) =:\mu ^{+}=1+C_{2}^{+}=1+\frac{\nu }{\alpha -1}>2.
\end{equation*}
gives the average contribution to $\Delta $ of the two components in the
special critical MC case. The $\Delta $s now constitute the interval lengths
of an alternating renewal process.

\subsubsection{\textbf{Time-reversal.}}

It is of importance to check whether or not detailed balance holds for the
MC under study here. Assume $X_{n}$ is recurrent. The catastrophe MC is not
time-reversible as detailed balance does not hold. Let $\overleftarrow{P}%
\neq P$ be the transition matrix of the process $X_{n}^{\leftarrow }$ which
is $X_{n}$ backward in time. With $^{\prime }$ denoting matrix transposition
and $D_{\mathbf{\pi }}=$diag$\left( \pi _{0},\pi _{1},...\right) $, we have 
\begin{equation*}
\overleftarrow{P}=D_{\mathbf{\pi }}^{-1}P^{\prime }D_{\mathbf{\pi }}.
\end{equation*}
It can be checked that the only non-null entries of $\overleftarrow{P}$ are
its first row with $\overleftarrow{P}_{0,0}=q_{0}$ and $\overleftarrow{P}%
_{0,x}=q_{x}\prod_{y=0}^{x-1}p_{y}$ if $x\geq 1$ and the lower diagonal
whose entries are all ones. Starting from $X_{0}^{\leftarrow }=x$, the
process $X_{n}^{\leftarrow }$ decays linearly till it first hits $0$ and
once in state $0$, $X_{n}^{\leftarrow }$ ends up jumping abruptly upward
(after some latency time if $q_{0}>0$). Equivalently, $X_{n}^{\leftarrow }$
undergoes a jump of amplitude $x\geq 0$ with probability $%
q_{x}\prod_{y=0}^{x-1}p_{y}$ before diminishing again and again to $0$. We
have; 
\begin{equation*}
\Bbb{P}\left( \tau _{0,0}-1=x\right) =q_{x}\prod_{y=0}^{x-1}p_{y}
\end{equation*}
so the jump's amplitude is the one of $\tau _{0,0}-1$ which is also the
height $H$ of an excursion. Note that $H=0$ with probability $q_{0}:$
whenever the walker started at $0$ remains at $0$ in the next step, this
event is considered as an excursion of length $1$ and height $0$.

The matrix $\overleftarrow{P}$ has the structure of an infinite-dimensional
Leslie (companion) matrix.

Clearly, in case of positive recurrence, $X_{n}^{\leftarrow }$ models the
forward recurrence time of the original process with $X_{n}^{\leftarrow }%
\overset{d}{\rightarrow }X_{\infty }^{\leftarrow }$ as $n\rightarrow \infty
. $ Note that $\mathbf{\pi }^{\prime }=\mathbf{\pi }^{\prime }\overleftarrow{%
P} $ ($\mathbf{\pi }$ is also the invariant measure for $X^{\leftarrow }$)
and so: $X_{\infty }\overset{d}{=}X_{\infty }^{\leftarrow }\overset{d}{\sim }%
\pi _{x}$, as required.

\subsubsection{\textbf{The scale (or harmonic) function.}}

In the recurrent case, the sample paths of $X_{n}$ are made of i.i.d.
excursions (the pieces of the sample paths between consecutive visits to $0$%
). The lengths of the excursions are $\tau _{0,0}.$ Let us look at their
heights $H$. As observed just before, $H\overset{d}{=}\tau _{0,0}-1,$
because $X$ grows linearly between consecutive visits to $0$ (if it grows).
Let us rapidly check this with the use of the scale function. The scale
function idea will appear useful.

Assume $X_{0}=x.$ Let $X_{n\wedge \tau _{x,0}}^{{}}$ stopping $X_{n}$ when
it first hits $0.$ Let us define the scale (or harmonic) function $\varphi $
of $X_{n}$ as the function which makes $Y_{n}\equiv \varphi \left(
X_{n\wedge \tau _{x,0}}^{{}}\right) $ a martingale. The function $\varphi $
is important because, as is well-known, for all $0<x<h,$ with $\tau _{x}$
the first hitting time of $\left\{ 0,h\right\} $ starting from $x$ (assuming 
$\varphi \left( 0\right) \equiv 0$)

\begin{equation*}
\Bbb{P}\left( X_{\tau _{x}}=h\right) =\Bbb{P}\left( \tau _{x,h}<\tau
_{x,0}\right) =\frac{\varphi \left( x\right) }{\varphi \left( h\right) }.
\end{equation*}
Using this remark, the event $H=h$ is realized when $\tau _{0,h}<\tau _{0,0}$
and $\tau _{h,h+1}>\tau _{h,0},$ the latter two events being independent.
Thus (recalling $\Bbb{P}\left( H=0\right) =q_{0}$): 
\begin{equation}
\Bbb{P}\left( H=h\right) =p_{0}\frac{\varphi \left( 1\right) }{\varphi
\left( h\right) }\left( 1-\frac{\varphi \left( h\right) }{\varphi \left(
h+1\right) }\right) \text{, }h\geq 1.  \label{h1}
\end{equation}
We clearly have $\sum_{h\geq 1}\Bbb{P}\left( H=h\right) =p_{0}$ because
partial sums are part of a telescoping series. But this is also $\Bbb{P}%
\left( H\geq h\right) =1/\varphi \left( h\right) $. It remains to compute $%
\varphi $ with $\varphi \left( 0\right) =0$. We wish to have: $\Bbb{E}%
_{x}\left( X_{n+1}\mid X_{n}=y\right) =y$, leading to 
\begin{equation*}
\varphi \left( x\right) =p_{x}\varphi \left( x+1\right) +q_{x}\varphi \left(
0\right) =p_{x}\varphi \left( x+1\right) ,\text{ }x\geq 1.
\end{equation*}
Thus, the searched `harmonic' (increasing) function is 
\begin{equation}
\varphi \left( x\right) =\frac{1}{\prod_{y=0}^{x-1}p_{y}}\text{, }x\geq 1,%
\text{ }\varphi \left( 0\right) \equiv 0.  \label{h2}
\end{equation}
Note $\varphi \left( 1\right) =1/p_{0}$ and $\varphi \left( x\right) $ is
diverging whenever the chain is recurrent. Equations (\ref{h1}) and (\ref{h2}%
) characterize the law of the excursion height of the random walker in the
recurrent case. Note 
\begin{equation*}
\Bbb{P}\left( H\geq h\right) =1/\varphi \left( h\right)
=\prod_{y=0}^{h-1}p_{y}=\Bbb{P}\left( \tau _{0,0}>h\right) ,
\end{equation*}
showing, as expected from the beginning, that $H\overset{d}{=}\tau _{0,0}-1.$

\textbf{Remark }(Doob transform)\textbf{:} Let $\overline{P}$ be obtained
from $P$ by removing its first and last column. With $\mathbf{\varphi :=}%
\left( \varphi \left( 1\right) ,\varphi \left( 2\right) ,...\right) $ and $%
D_{\mathbf{\varphi }}=$diag$\left( \mathbf{\varphi }\right) $, the Markov
chain with stochastic transition matrix:

\begin{equation*}
\overline{P}_{\mathbf{\varphi }}=D_{\mathbf{\varphi }}^{-1}\overline{P}D_{%
\mathbf{\varphi }}
\end{equation*}
is $X$ (with transition matrix $P$) conditioned to hit first $\infty $
before $0$. In our context: 
\begin{equation*}
\overline{P}_{\mathbf{\varphi }}=\left[ 
\begin{array}{llllll}
0 & 1 &  &  &  & \cdots \\ 
0 & 0 & 1 &  &  & \cdots \\ 
\vdots & \vdots & \ddots & \ddots &  & \cdots \\ 
0 & 0 & \cdots & \ddots & 1 & \cdots \\ 
\vdots & 0 & \cdots &  & 0 & \ddots \\ 
\vdots & \vdots &  &  &  & \ddots
\end{array}
\right] .\text{ }\Box
\end{equation*}

\subsubsection{\textbf{Probability of extinction.}}

Consider now the same Markov chain but assume now that $p_{0}=0$, $q_{0}=1.$
In this case, the state $0$ is absorbing. Consider then the restriction $%
\overline{P}$ of matrix $P$ to the states $\left\{ 1,2,...\right\} $. Let $%
\phi _{x}$, $x\geq 1$ be the probabilities that state $0$ is hit in finite
time given the chain started originally at $x$. Let $\mathbf{\phi }\equiv
\left( \phi _{1},\phi _{2},..\right) ^{\prime }$ be the column-vector of
these absorption probabilities. Let $\mathbf{q}\equiv \left(
q_{1},q_{2},..\right) ^{\prime }.$ Then $\mathbf{\phi }$ is the smallest
non-negative solution to $\mathbf{\phi =q}+\overline{P}\mathbf{\phi }$ whose
formal solution is: $\mathbf{\phi }=\left( I-\overline{P}\right) ^{-1}%
\mathbf{q.}$ The $\phi _{x}$s obeys the recurrence: $\phi
_{x}=q_{x}+p_{x}\phi _{x+1},$ else $1-\phi _{x+1}=\frac{1}{p_{x}}\left(
1-\phi _{x}\right) $. All $\phi _{x}$ can therefore be expressed in terms of 
$\phi _{1}$, leading simply to: 
\begin{equation*}
1-\phi _{x}=\frac{1}{\prod_{y=1}^{x-1}p_{y}}\left( 1-\phi _{1}\right) .
\end{equation*}
The formal solution is also $\mathbf{\phi }=\left( I-\overline{P}\right)
^{-1}\mathbf{q}$, involving the resolvent of $\overline{P}$. Because $\left(
I-\overline{P}\right) ^{-1}$is computable with upper triangular structure, 
\begin{equation*}
\left( I-\overline{P}\right) ^{-1}=\left[ 
\begin{array}{llllll}
1 & p_{1} & p_{1}p_{2} & p_{1}p_{2}p_{3} & p_{1}p_{2}p_{3}p_{4} & \cdots \\ 
& 1 & p_{2} & p_{2}p_{3} & p_{2}p_{3}p_{4} & \cdots \\ 
&  & 1 & p_{3} & p_{3}p_{4} & \cdots \\ 
&  &  & 1 & p_{4} & \cdots \\ 
&  &  &  & 1 & \cdots \\ 
&  &  &  &  & \ddots
\end{array}
\right] ,
\end{equation*}
$\phi _{x}$ thus takes the alternative form 
\begin{equation*}
\phi _{x}=\sum_{y\geq x}q_{y}\prod_{y^{\prime }=x}^{y-1}p_{y^{\prime }}.
\end{equation*}
If $C_{1}\equiv \sum_{y\geq 1}q_{y}=\infty $, then $u_{x}=%
\prod_{y=1}^{x-1}p_{y}\rightarrow 0:$ the restriction $\phi _{x}\in \left[
0,1\right] $ forces $\phi _{1}=1$ and so $\phi _{x}=1$ for all $x<\infty $:
The state $0$ is hit with probability $1$, starting from $x,$ for all $%
x<\infty .$ The Markov chain is recurrent.

But, if $C_{1}<\infty $, then we can take $\phi _{1}<1$ so long as $\phi
_{x}\geq 0$ for all $x\geq 1.$ The minimal solution occurs when $1-\phi _{1}=%
\Bbb{P}\left( \tau _{1,0}=\infty \right) =\prod_{y\geq 1}p_{y}>0,$ leading
to the alternative expression of $\phi _{x}$: 
\begin{equation}
\phi _{x}=1-\frac{\prod_{y\geq 1}p_{y}}{\prod_{y=1}^{x-1}p_{y}}%
=1-\prod_{y\geq x}p_{y}.  \label{3}
\end{equation}
In this case, $\phi _{x}<1$ for $x\geq 1$ and the absorbed random walker
started at $x$ avoids $0$ with positive probability (a transience case for
the original reflected Markov chain). Note that $x<x^{\prime }\Rightarrow
\phi _{x}>\phi _{x^{\prime }}.$ We can extend Eq. (\ref{3}) to $x=0$,
because, by first-step analysis, 
\begin{equation*}
\phi _{0}:=\Bbb{P}\left( \tau _{0,0}<\infty \right) =\phi _{0,0}\left(
1\right) =p_{0}\Bbb{P}\left( \tau _{1,0}<\infty \right) +q_{0}=p_{0}\phi
_{1}+q_{0},
\end{equation*}
leading to 
\begin{equation*}
\phi _{0}=1-p_{0}\left( 1-\phi _{1}\right) =1-\prod_{y\geq 0}p_{y}.
\end{equation*}
\newline

To summarize, we have:

$\left( i\right) $ If $C_{1}<\infty ,$ the MC is transient and, with 
\begin{equation*}
\tau _{x,0}=\inf \left( n\geq 1:X_{n}=0\mid X_{0}=x\right) ,\text{ }x\geq 0,
\end{equation*}
$\Bbb{P}\left( \tau _{x,0}=\infty \right) =\prod_{y\geq x}p_{y}>0$. The
chain $X$ started at $x\geq 0$ has probability $\phi _{x}=1-\prod_{y\geq
x}p_{y}<1$ to undergo a first extinction.

$\left( ii\right) $ If $C_{1}=\infty ,$ the MC is recurrent with $\Bbb{P}%
\left( \tau _{x,0}=\infty \right) =0$. Moreover, it is: null recurrent if $%
C_{2}=\infty $, positive recurrent if $C_{2}<\infty .$

Due to irreducibility (because $p_{x}$ and $q_{x}>0$, for all $x\geq 1$),
states are either all transient or recurrent.\newline

APPLICATION: for the model (\ref{0a}), 
\begin{equation*}
\Bbb{P}\left( \tau _{x,0}=\infty \right) >0\Leftrightarrow \beta >1.
\end{equation*}

\subsubsection{\textbf{Times to collapse (first extinction).}}

How long does it take, starting from $x\geq 1$, to first hit $0?$ We give
here some insight on the way to compute the law of these first times to
collapse. With $x\geq 1$, let thus $\tau _{x,0}$ be the time it takes to
first hit $0$, starting from $X_{0}=x\geq 1$. With $\tau _{x+1,0}^{\prime }$
a statistical copy of $\tau _{x+1,0}$, from first-step analysis, we clearly
have: 
\begin{equation*}
\tau _{x,0}\overset{d}{=}\left( 1-B_{x}\right) \cdot 1+B_{x}\cdot \left(
1+\tau _{x+1,0}^{\prime }\right) ,
\end{equation*}
where $B_{x}$ is a Bernoulli random variable with $\Bbb{P}\left(
B_{x}=1\right) =p_{x}.$ Therefore with $\phi _{x,0}\left( z\right) =\Bbb{E}%
\left( z^{\tau _{x,0}}\right) ,$ $\phi _{x}\left( z\right) $ obeys the
recurrence $\phi _{x,0}\left( z\right) =q_{x}z+p_{x}z\phi _{x+1,0}\left(
z\right) $, with initial condition $\phi _{1,0}\left( z\right) $: again, $%
\phi _{x,0}\left( z\right) $ can easily be deduced once $\phi _{1,0}\left(
z\right) $ is known. The recurrence is also 
\begin{equation*}
1-\phi _{x+1,0}\left( z\right) =\frac{z-1}{p_{x}z}+\frac{1-\phi _{x,0}\left(
z\right) }{p_{x}z}.
\end{equation*}
The full pgf of $\tau _{x,0}$ follows by recurrence. When $z=1$, $\phi
_{x,0}\left( 1\right) =\phi _{x}=\Bbb{P}\left( \tau _{x,0}<\infty \right) $
are the absorption probabilities already computed.

With $\mathbf{\phi }\left( z\right) =\left( \phi _{1,0}\left( z\right) ,\phi
_{2,0}\left( z\right) ,...\right) ^{\prime }$ the column-vector of the $\phi
_{x,0}\left( z\right) $, and $\mathbf{q}=\left( q_{1},q_{2},...\right)
^{\prime }$ the column-vector of the $q_{x},$ $\mathbf{\phi }\left( z\right) 
$ solves: 
\begin{equation}
\mathbf{\phi }\left( z\right) =z\mathbf{q}+z\overline{P}\mathbf{\phi }\left(
z\right) ,  \label{4}
\end{equation}
whose formal solution is $\mathbf{\phi }\left( z\right) =z\left( I-z%
\overline{P}\right) ^{-1}\mathbf{q}$, involving the resolvent of $\overline{P%
}$ which is 
\begin{equation*}
\left( I-z\overline{P}\right) ^{-1}=\left[ 
\begin{array}{llllll}
1 & zp_{1} & z^{2}p_{1}p_{2} & z^{3}p_{1}p_{2}p_{3} & 
z^{4}p_{1}p_{2}p_{3}p_{4} & \cdots \\ 
& 1 & zp_{2} & z^{2}p_{2}p_{3} & z^{3}p_{2}p_{3}p_{4} & \cdots \\ 
&  & 1 & zp_{3} & z^{2}p_{3}p_{4} & \cdots \\ 
&  &  & 1 & zp_{4} & \cdots \\ 
&  &  &  & 1 & \cdots \\ 
&  &  &  &  & \ddots
\end{array}
\right] .
\end{equation*}
We get 
\begin{equation*}
\phi _{x,0}\left( z\right) =\sum_{y\geq x}q_{y}z^{y-x+1}\prod_{y^{\prime
}=x}^{y-1}p_{y^{\prime }}.
\end{equation*}
Equivalently, with $x\geq 1$, the pmf of $\tau _{x,0}$ reads 
\begin{equation*}
\Bbb{P}\left( \tau _{x,0}=k\right) =q_{k+x-1}\prod_{y^{\prime
}=x}^{k+x-2}p_{y^{\prime }}\text{, }k\geq 1.
\end{equation*}
The above first-step analysis clearly also makes sense if one starts from $%
x=0$ and the recurrence $\phi _{x,0}\left( z\right) =q_{x}z+p_{x}z\phi
_{x+1,0}\left( z\right) $ also holds when $x=0.$ We conclude 
\begin{equation*}
\phi _{0,0}\left( z\right) =q_{0}z+\sum_{y\geq
1}q_{y}z^{y+1}\prod_{y^{\prime }=0}^{y-1}p_{y^{\prime }}=\sum_{y\geq
0}q_{y}z^{y+1}\prod_{y^{\prime }=0}^{y-1}p_{y^{\prime }},
\end{equation*}
consistently with (\ref{1a}).\newline

APPLICATION: for the model (\ref{0a}), 
\begin{equation*}
\Bbb{E}\left( \tau _{x,0}\right) =\infty \Leftrightarrow \beta >1\text{
(transience) or }\beta =1\text{ and }\alpha \leq 1\text{ (null recurrence).}
\end{equation*}
If $\beta <1$ or $\beta =1$ and $\alpha >1$ (positive recurrence), $\Bbb{E}%
\left( \tau _{x,0}\right) =\phi _{x,0}^{\prime }\left( 1\right) <\infty $.

\subsubsection{\textbf{Transience versus recurrence.}}

We here discuss the criterion for recurrence or transience of the general
catastrophe Markov chain.

- When $C_{1}=\infty ,$ the recurrent chain started at $x$ first hits $0$
with probability $1$ and returns infinitely often to $0$. Given $X_{0}=x$,
with $\mathcal{N}_{x,y}\equiv \sum_{n\geq 0}\mathbf{1}\left( X_{n}=y\right)
, $ the number of visits to state $y$, then $\mathcal{N}_{x,y}=\infty ,$ $%
\Bbb{P}_{x}-$almost surely. If $\tau _{x,x}$ is the first return time at $x$%
, then $\Bbb{P}\left( \tau _{x,x}<\infty \right) =1.$ Furthermore, with $%
\mathcal{N}_{x,y}\equiv \sum_{n=0}^{\tau _{x,x}}\mathbf{1}\left(
X_{n}=y\right) $ the number of visits to state $y$ before the first return
time to state $x$, then: $\Bbb{E}\left( \mathcal{N}_{x,y}\right) =\frac{\pi
_{y}}{\pi _{x}}$ and by the Chacon-Ornstein limit ratio ergodic theorem: 
\begin{equation*}
\frac{\sum_{n=0}^{N}\mathbf{1}\left( X_{n}=y\right) }{\sum_{n=0}^{N}\mathbf{1%
}\left( X_{n}=x\right) }\underset{N\nearrow \infty }{\rightarrow }\frac{\pi
_{y}}{\pi _{x}},\text{ }\Bbb{P}_{x}-\text{almost surely.}
\end{equation*}
Starting in particular from $x=0$, a recurrent chain is made of infinitely
many independent and identically distributed (i.i.d.) excursions which are
the sample paths of $\left( X_{n};n\geq 0\right) $ between consecutive
visits to state $0$. We have: $\Bbb{E}\left( \mathcal{N}_{0,x}\right) =\frac{%
\pi _{x}}{\pi _{0}}=\prod_{y=0}^{x-1}p_{y}.$ When the chain is positive
recurrent ($C_{2}<\infty $) the expected time elapsed between consecutive
visits to $0$ is finite and equal to $\Bbb{E}\left( \tau _{0,0}\right)
\equiv \mu =1/\pi _{0}=1+C_{2}$, whereas this expected time is infinite when
the chain is null recurrent.

- When $C_{1}<\infty ,$ the state $x\geq 0$ is transient.

Thus, $\mathcal{N}_{x,x}<\infty ,$ $\Bbb{P}_{x}-$almost surely and $\Bbb{P}%
\left( \mathcal{N}_{x,x}=k\right) =\left( 1-\alpha _{x}\right) \alpha
_{x}^{k-1}$ where $\alpha _{x}=\Bbb{P}\left( \tau _{x,x}<\infty \right) <1.$

Also, starting from $x\geq 1$, the walker hits at once $0$ with probability $%
\phi _{x}=$ $\Bbb{P}\left( \tau _{x,0}<\infty \right) $ and given this
occurred it undergoes a number $\mathcal{N}_{e}=k$ excursions with $\Bbb{P}%
\left( \mathcal{N}_{e}=k\right) =\left( 1-\phi _{0,0}\right) \phi _{0,0}^{k}$%
, $k\geq 0$ before drifting to $\infty $ for ever. With probability $1-\phi
_{x}$, the walker drifts to $\infty $ without ever visiting $0$ and $%
\mathcal{N}_{e}=0$. The time $\tau _{x}^{\left( d\right) }$ at which the
walker (initially at $x$) starts drifting linearly to $\infty $ for ever is
thus:

\begin{equation*}
\tau _{x}^{\left( d\right) }:=\left( \tau _{x,0}+\sum_{k=1}^{\mathcal{N}%
_{e}}\tau _{0,0}^{\left( k\right) }\right) \mathbf{1}\left( \tau
_{x,0}<\infty \right)
\end{equation*}
where the $\tau _{0,0}^{\left( k\right) }$s are iid copies of $\tau _{0,0}.$

\subsubsection{\textbf{First passage times and Green kernel.}}

Let $\tau _{x,y}$ be the first passage time at $y\neq x$ when the process $X$
is started at $x$. We wish here to briefly derive an exact formal formula
for the law of $\tau _{x,y}$, making use of the Green function of a MC. Let 
\begin{equation*}
\phi _{x,y}\left( z\right) \equiv \sum_{k=1}^{\infty }z^{k}\Bbb{P}\left(
\tau _{x,y}=k\right)
\end{equation*}
be the generating function of the law of $\tau _{x,y}.$ Then, with 
\begin{equation*}
g_{x,y}\left( z\right) \equiv \sum_{n=0}^{\infty }z^{n}\Bbb{P}_{x}\left(
X_{n}=y\right) =\sum_{n=0}^{\infty }z^{n}P^{n}\left( x,y\right) =\left(
I-zP\right) ^{-1}\left( x,y\right)
\end{equation*}
the generating function of $P^{n}\left( x,y\right) $ (the Green potential
function of the chain), using $P^{n}\left( x,y\right) =\sum_{m=0}^{n}\Bbb{P}%
\left( \tau _{x,y}=m\right) P^{n-m}\left( y,y\right) $, we easily get the
expression: 
\begin{equation*}
\phi _{x,y}\left( z\right) =\frac{g_{x,y}\left( z\right) }{g_{y,y}\left(
z\right) }.
\end{equation*}
In particular, 
\begin{equation*}
\phi _{x,0}\left( z\right) =\frac{g_{x,0}\left( z\right) }{g_{0,0}\left(
z\right) }\text{ and }\phi _{0,x}\left( z\right) =\frac{g_{0,x}\left(
z\right) }{g_{x,x}\left( z\right) }
\end{equation*}
are the generating functions of $\tau _{x,0}$ and $\tau _{0,x}.$

The pgf $\phi _{x,x}\left( z\right) $ of the first-return time $\tau _{x,x}$
to state $x$ satisfies 
\begin{equation*}
\phi _{x,x}\left( z\right) =\frac{g_{x,x}\left( z\right) -1}{g_{x,x}\left(
z\right) }=1-\frac{1}{g_{x,x}\left( z\right) }
\end{equation*}
where $g_{x,x}\left( z\right) =\sum_{n=0}^{\infty }z^{n}\Bbb{P}_{x}\left(
X_{n}=x\right) =\sum_{n=0}^{\infty }z^{n}P^{n}\left( x,x\right) $ is the
Green function at $x.$ Recall (in the positive recurrent case with $\phi
_{x,x}^{\prime }\left( 1\right) <\infty $) 
\begin{equation*}
\Bbb{E}\left( \tau _{x,x}\right) =\phi _{x,x}^{\prime }\left( 1\right)
=1/\pi _{x}.
\end{equation*}

In the specific MCC the resolvent can easily be computed. We find:

\begin{equation*}
g_{x,x}\left( z\right) =\left( I-zP\right) ^{-1}\left( x,x\right) =\frac{%
1-\sum_{x^{\prime }=0}^{x-1}z^{x^{\prime }+1}q_{x^{\prime
}}\prod_{y=1}^{x^{\prime }-1}p_{y}}{1-\sum_{x^{\prime }\geq 0}z^{x^{\prime
}+1}q_{x^{\prime }}\prod_{y=0}^{x^{\prime }-1}p_{y}}
\end{equation*}
giving the Green kernel as 
\begin{eqnarray*}
g_{x,y}\left( z\right) &=&g_{x,x}\left( z\right) z^{y-x}\prod_{y^{\prime
}=x}^{y-1}p_{y^{\prime }}\text{ if }y\geq x \\
g_{x,y}\left( z\right) &=&\left( g_{x,x}\left( z\right) -1\right)
z^{y-x}/\prod_{y^{\prime }=y}^{x-1}p_{y^{\prime }}\text{ if }0\leq y<x.
\end{eqnarray*}
Note that the denominator term of $g_{x,x}\left( z\right) $ is $1-\phi
_{0,0}\left( z\right) $ where $\phi _{0,0}\left( z\right) $ is the pgf of $%
\tau _{0,0}$ computed above in the special case. This yields 
\begin{equation*}
\phi _{x,x}\left( z\right) =\frac{\sum_{x^{\prime }\geq x}z^{x^{\prime
}+1}q_{x^{\prime }}\prod_{y=0}^{x^{\prime }-1}p_{y}}{1-\sum_{x^{\prime
}=0}^{x-1}z^{x^{\prime }+1}q_{x^{\prime }}\prod_{y=0}^{x^{\prime }-1}p_{y}}.
\end{equation*}
In particular, 
\begin{equation*}
\alpha _{x}=\Bbb{P}\left( \tau _{x,x}<\infty \right) =\phi _{x,x}\left(
1\right) .
\end{equation*}
Note also that 
\begin{equation*}
g_{0,0}\left( z\right) =1/\left( 1-\phi _{0,0}\left( z\right) \right)
=\sum_{n=0}^{\infty }z^{n}\Bbb{P}_{0}\left( X_{n}=0\right) ,
\end{equation*}
where $\Bbb{P}_{0}\left( X_{n}=0\right) $ is the contact probability at $0$
at time $n$.\newline

APPLICATION (contact probability at $0$): With $\phi _{0,0}\left( z\right)
=z\left( q_{0}+p_{0}\psi _{0}\left( z\right) \right) $ in the critical
special case.

- If $0<\alpha <1$, using $\psi _{0}\left( z\right) \underset{z\downarrow 1}{%
\sim }1-\frac{\Gamma \left( \nu +1\right) \Gamma \left( 1+\alpha \right) }{%
\Gamma \left( \nu +1-\alpha \right) }\left( 1-z\right) ^{\alpha },$ we have 
\begin{equation*}
1-\phi _{0,0}\left( z\right) \underset{z\downarrow 1}{\sim }p_{0}\left(
1-\psi _{0}\left( z\right) \right) \sim p_{0}\frac{\Gamma \left( \nu
+1\right) \Gamma \left( 1-\alpha \right) }{\Gamma \left( \nu +1-\alpha
\right) }\left( 1-z\right) ^{\alpha }
\end{equation*}
showing by singularity analysis that, in the null recurrent case (algebraic
decay of the contact probability) 
\begin{equation*}
\Bbb{P}_{0}\left( X_{n}=0\right) \underset{n\uparrow \infty }{\sim }\frac{%
\Gamma \left( \nu +1-\alpha \right) }{p_{0}\Gamma \left( \nu +1\right)
\Gamma \left( 1-\alpha \right) \Gamma \left( \alpha \right) }n^{-\left(
1-\alpha \right) }.
\end{equation*}
When $\alpha =1-\varepsilon $ ($\varepsilon >0$ small), the constant in
front of $n^{-\left( 1-\alpha \right) }$ vanishes like $\varepsilon /\left(
p_{0}\nu \right) $.

- When $\alpha =1$, using $\psi _{0}\left( z\right) \underset{z\downarrow 1}{%
\sim }1-\left( 1-z\right) -\nu \left( 1-z\right) \log \left( 1-z\right) $%
\begin{equation*}
1-\phi _{0,0}\left( z\right) =1-z\left( q_{0}+p_{0}\psi _{0}\left( z\right)
\right) \underset{z\downarrow 1}{\sim }\left( 1-z\right) \left(
1+p_{0}\right) +p_{0}\nu \left( 1-z\right) \log \left( 1-z\right) ,
\end{equation*}
with logarithmic singularity, showing by singularity analysis, that 
\begin{equation*}
\Bbb{P}_{0}\left( X_{n}=0\right) \underset{n\uparrow \infty }{\sim }\frac{1}{%
p_{0}\nu \log n}.
\end{equation*}

- When $\alpha >1$ (positive recurrence case), using $\psi _{0}\left(
z\right) \underset{z\downarrow 1}{\sim }1-\frac{\nu }{\alpha -1}\left(
1-z\right) +O\left( 1-z\right) ^{\alpha },$ 
\begin{equation*}
1-\phi _{0,0}\left( z\right) \underset{z\downarrow 1}{\sim }1-z\left(
q_{0}+p_{0}\left( 1-\frac{\nu }{\alpha -1}\left( 1-z\right) \right) \right)
\sim \left( 1-z\right) \left( 1+p_{0}\frac{\nu }{\alpha -1}\right)
\end{equation*}
showing, as required that $\Bbb{P}_{0}\left( X_{n}=0\right) \underset{%
n\uparrow \infty }{\rightarrow }\pi _{0}=1/\left( 1+p_{0}\frac{\nu }{\alpha
-1}\right) .$ When $\alpha =1+\varepsilon $ ($\varepsilon >0$ small): $\pi
_{0}\sim \varepsilon /\left( p_{0}\nu \right) $, just like when $\alpha <1.$%
\newline

\textbf{Remark }(spectral aspects of $P$): for the choice of $\left(
p_{x},q_{x}\right) $ as from (\ref{0a}), in the critical case $\beta =1$,
the operator $P$ can easily be checked not to be compact, nor even
quasi-compact, \cite{Gos}. Looking for sequences $\mathbf{u}^{\prime
}=\left( u_{0},u_{1},u_{2},...\right) $\ solving 
\begin{equation*}
\left( I-zP\right) \mathbf{u=0}
\end{equation*}
yields 
\begin{equation*}
u_{x}=\frac{u_{0}}{z^{x}\prod_{y=0}^{x-1}p_{y}}\left( 1-\sum_{x^{\prime
}=0}^{x-1}z^{x^{\prime }+1}q_{x^{\prime }}\prod_{y=0}^{x^{\prime
}-1}p_{y}\right) ,\text{ }x\geq 0
\end{equation*}
which are well-defined in the Banach space $c_{0}=\left\{ \mathbf{u}%
:u_{x}\rightarrow 0\emph{\ }\text{as}\emph{\ }x\rightarrow \infty \right\} $%
\ whenever $\func{Re}\left( z\right) \geq 1.$\ We conclude that for each $z$
obeying $\func{Re}\left( z\right) \geq 1$, there are solutions $\mathbf{u}$\
(eigen-states), vanishing at $\infty $, to $\left( I-zP\right) \mathbf{u}=0,$%
\ defined up to a multiplicative constant.\ Letting $\lambda =z^{-1},$\ we
get that the closed disk of $\Bbb{C}$\ centered at $\left( 1/2,0\right) $\
with radius $1/2$\ (which is: $\func{Re}\left( \lambda ^{-1}\right) \geq 1$)
constitutes the point spectrum of $P.$\ When $\lambda $\ belongs to the
latter disk with radius $1/2$, $\left( \lambda I-P\right) ^{-1}$\ does not
exist. Because there is a continuum of eigenvalues in the latter disk, the
corresponding operator $P$\ has no spectral gap. The points $\lambda $\
belonging to the complementary of the latter disk in the unit disk centered
at $0$\ constitute the continuous spectrum where $\left( \lambda I-P\right)
^{-1}$\ exists but is not defined on the whole space $c_{0}.$ The $\lambda $%
s belonging to the open domain $\left| \lambda ^{-1}\right| <1$ (the
complementary of the unit disk of $\Bbb{C}$ centered at $0$) are regular
points of $P$ for which $\left( \lambda I-P\right) ^{-1}$ exists, is defined
on the whole space $c_{0}$ and is bounded. $\Box $

\subsection{\textbf{Including rates and time change: continuous-time}}

Consider now the continuous-time catastrophe MC $\overline{X}\left( t\right)
,$ so with $\overline{X}\left( t\right) $ generated by the transition rate
matrix: $Q=D_{\mathbf{r}}\left( P-I\right) $ defined in (\ref{0b}). We now
specify the results of the latter Section to this special MC in
continuous-time, including the opportunity of a holding rate in the
transition matrix of the DT catastrophe process. This special MC deserves
special interest in particular because it is, to a large extent, amenable to
exact analytic computations. The study of $\overline{X}$ should take into
account the three parameters $\alpha ,\beta >0,\nu >-1$, $\alpha <\nu +1$,
together with the rates\textbf{\ }$r_{x}=r_{0}\left( x+1\right) ^{\lambda
}>0 $, $x\geq 0$.

\subsubsection{\textbf{Invariant measure.}}

We now investigate the way the invariant measure is modified by the
adjunction of holding rates. Let $D_{\mathbf{r}}=$diag$\left(
r_{0},r_{1},r_{2},...\right) .$ We have $Q=D_{\mathbf{r}}\left( P-I\right) .$
Let $\overline{\mathbf{\pi }}$ be the invariant measure associated to $Q$,
when it exists. It should solve $\mathbf{0}^{\prime }=\overline{\mathbf{\pi }%
}^{\prime }Q,$ and we get: 
\begin{equation*}
r_{0}p_{0}\overline{\pi }_{0}=\sum_{x\geq 1}\overline{\pi }_{x}q_{x}r_{x}%
\text{ and }\overline{\pi }_{x}=\overline{\pi }_{0}r_{0}\frac{%
\prod_{y=0}^{x-1}p_{y}}{r_{x}}\text{, }x\geq 1.
\end{equation*}
Using the second equation, the first equation is satisfied whenever 
\begin{equation*}
\sum_{x\geq 1}q_{x}r_{x}\frac{\prod_{y=1}^{x-1}p_{y}}{r_{x}}=\sum_{x\geq
1}\left( u_{x}-u_{x+1}\right) =1,
\end{equation*}
so again when $u_{x}=\prod_{y=1}^{x-1}p_{y}\underset{x\rightarrow \infty }{%
\rightarrow }0$ which is fulfilled if and only if $C_{1}\equiv
\sum_{y=1}^{\infty }q_{y}=\infty $ (the recurrence condition for $X$): the
time change leading from $X$ to $\overline{X}$ does not change the road map
of $X$ so the recurrence criterion is identical for both $X$ and $\overline{X%
}.$ However, the criterion for positive recurrence is modified.

Indeed, if in addition $\overline{C}_{2}\equiv \sum_{x\geq 1}\left( \frac{%
r_{x}}{r_{0}}\right) ^{-1}\prod_{y=0}^{x-1}p_{y}<\infty $, then $\overline{%
\pi }_{0}=\frac{1}{1+\overline{C}_{2}}\in \left( 0,1\right) $ and then the
MC is positive recurrent with invariant probability measure 
\begin{equation*}
\overline{\pi }_{x}=\overline{\pi }_{0}r_{0}\frac{\prod_{y=0}^{x-1}p_{y}}{%
r_{x}}\text{, }x\geq 0.
\end{equation*}
Else if $\overline{C}_{2}=\infty $, the MC $\overline{X}$ is null recurrent.

APPLICATION: When dealing with the special continuous-time Markov chain
(CTMC), we conclude:

- $C_{1}<\infty $ if and only if $\beta >1$: in this case the MC is
transient with no invariant measure.

- If $\beta <1,$ then $\overline{C}_{2}<\infty :$ the MC is positive
recurrent. Furthermore, for large $x$: 
\begin{equation*}
\overline{\pi }_{x}\sim x^{-\lambda }\prod_{y=0}^{x-1}\left( 1-\frac{\alpha 
}{\nu +y^{\beta }}\right) \sim x^{-\lambda }e^{-\alpha \sum^{x}y^{-\beta
}}\sim x^{-\lambda }e^{-\frac{\alpha }{1-\beta }x^{1-\beta }}
\end{equation*}
with stretched exponential behaviour.

- In the critical case $\beta =1$: 
\begin{equation*}
\overline{\pi }_{x}\sim x^{-\lambda }\prod_{y=0}^{x-1}\left( 1-\frac{\alpha 
}{\nu +y}\right) \sim x^{-\left( \lambda +\alpha \right) }
\end{equation*}
with power-law$\left( \alpha +\lambda \right) $ behaviour. The chain is
recurrent. From this one can guess that the CT chain is positive recurrent
if $\alpha +\lambda >1$, null recurrent if $\alpha +\lambda \leq 1.$ In
particular, the chain is always positive recurrent if $\lambda >1.$

\subsubsection{\textbf{Scale function.}}

$\overline{X}\left( t\right) $ is again the size at time $t$ of the current
population before the last catastrophic event. In the recurrent case ($%
C_{1}=\infty $), $\overline{X}$ is again made of i.i.d. excursions, whose
height $\overline{H}$ has the same law as the one for $X:$ indeed, one can
easily check that the scale function $\overline{\varphi }$ of $\overline{X}$
solving $Q\overline{\varphi }=0$, $\overline{\varphi }_{0}=0,$ coincide with
the scale function $\varphi $ of $X$ solving $P\varphi =\varphi $, $\varphi
_{0}=0.$ Because the scale function determines the height's law, we have the
claimed statement that the height's law is left unchanged by time
substitution.

\subsubsection{\textbf{Excursions lengths.}}

However, because of the time change, the times $\overline{\tau }_{0,0}$
between consecutive visits to $0$ (the excursions lengths of $\overline{X}$)
are of course very different from $\tau _{0,0}$ (statistically longer). With 
$H\overset{d}{=}\tau _{0,0}-1$ and $E_{y}\sim $Exp$\left( r_{y}\right) $%
\begin{equation*}
\overline{\tau }_{0,0}=\sum_{y=0}^{H}E_{y}
\end{equation*}
with the $E_{x}$s mutually independent and independent of $H$. [\cite{And},
page $21$].

Given the height of an excursion is $h$, with 
\begin{equation*}
c_{h,x}=\prod_{y=0:y\neq x}^{h}\frac{r_{y}}{r_{y}-r_{x}}
\end{equation*}
we have 
\begin{equation*}
\Bbb{P}\left( \overline{\tau }_{0,0}>t\mid H=h\right) =\Bbb{P}\left(
\sum_{x=0}^{h}E_{x}>t\right) =\sum_{x=0}^{h}c_{h,x}e^{-r_{x}t}.
\end{equation*}
and with $x_{*}=\arg \min \left( r_{x}:x=0,...,h\right) :$ $\Bbb{P}\left( 
\overline{\tau }_{0,0}>t\mid H=h\right) \sim c_{h,x_{*}}e^{-r_{x_{*}}t}$ as $%
t$ is large. Note $x_{*}=0$ if $\lambda >0$, $=h$ if $\lambda <0.$

Consider the critical case ($\beta =1$) with $\lambda <1$. Estimating the
Laplace-transform of $\overline{\tau }_{0,0}$ yields 
\begin{equation*}
\Bbb{E}\left( e^{-p\overline{\tau }_{0,0}}\right) =\sum_{h\geq 0}\Bbb{P}%
\left( H=h\right) \Bbb{E}\left( e^{-p\sum_{x=0}^{h}E_{x}}\right)
=\sum_{h\geq 0}\Bbb{P}\left( H=h\right) \prod_{x=0}^{h}\frac{1}{1+p/r_{x}}
\end{equation*}

\begin{equation*}
\underset{p\rightarrow 0}{\sim }\sum_{h\geq 0}\Bbb{P}\left( H=h\right)
e^{-p\sum_{x=0}^{h}1/r_{x}}\sim \sum_{h\geq 0}\Bbb{P}\left( H=h\right) e^{-%
\frac{p}{r_{0}\left( 1-\lambda \right) }h^{1-\lambda }}=\Bbb{E}\left( e^{-%
\frac{p}{r_{0}\left( 1-\lambda \right) }H^{1-\lambda }}\right)
\end{equation*}
suggesting $\overline{\tau }_{0,0}$ is tail-equivalent to $\frac{1}{%
r_{0}\left( 1-\lambda \right) }H^{1-\lambda }$. Recalling $\Bbb{P}\left(
H\geq h\right) \approx h^{-\alpha }$, we get 
\begin{equation*}
\Bbb{P}\left( \overline{\tau }_{0,0}>t\right) \approx t^{-\alpha /\left(
1-\lambda \right) }.
\end{equation*}
Clearly, $\Bbb{E}\left( \overline{\tau }_{0,0}\right) <\infty $ if and only
if $\alpha /\left( 1-\lambda \right) >1$ which, when $\lambda <1$, is the
positive-recurrence criterion.

If $\lambda >1,$\ 
\begin{equation*}
\sum_{x=0}^{h}1/r_{x}\sim \frac{1}{r_{0}\left( \lambda -1\right) }\left(
1-h^{-\left( \lambda -1\right) }\right) \nearrow \frac{1}{r_{0}\left(
\lambda -1\right) }
\end{equation*}
suggesting $\overline{\tau }_{0,0}\sim \frac{1}{r_{0}\left( \lambda
-1\right) }\left( 1-H^{-\left( \lambda -1\right) }\right) $\ with finite
mean, always. With $\mu _{h}=\sum_{x=0}^{h}1/r_{x}$ ($\nearrow \mu =\frac{1}{%
r_{0}\left( \lambda -1\right) }$ as $h\rightarrow \infty $), this is
consistent with [\cite{Jan}, Theorem $5.1,$ $(ii)$] 
\begin{equation*}
\Bbb{P}\left( \overline{\tau }_{0,0}>t\mid H=h\right) =\Bbb{P}\left(
\sum_{x=0}^{h}E_{x}\geq t\right) \leq ee^{-t/\mu _{h}}<ee^{-t/\mu }\text{
when }t\text{ is large}.
\end{equation*}

\subsubsection{\textbf{Excursions heights.}}

As observed just before, they are the same as the ones of the discrete-time
MC (DTMC).

\subsubsection{\textbf{Time to first extinction.}}

In the continuous-time setting, with $\left( E_{y},\text{ }y\geq x\right) $
mutually independent rvs with distribution Exp$\left( r_{y}\right) ,$
independent of $\tau _{x,0}$, we clearly have 
\begin{equation*}
\overline{\tau }_{x,0}=\sum_{y=x}^{\tau _{x,0}+x-1}E_{y},
\end{equation*}
consistently with the first-step analysis giving: 
\begin{equation*}
\overline{\tau }_{x,0}\overset{d}{=}\left( 1-B_{x}\right) \cdot
E_{x}+B_{x}\cdot \left( E_{x}+\overline{\tau }_{x+1,0}^{\prime }\right) .
\end{equation*}

\subsubsection{\textbf{Time }$\overline{\tau }_{x}^{\left( d\right) }$%
\textbf{\ at which the walker (initially at }$x$\textbf{) starts drifting
linearly to }$\infty $\textbf{\ for ever.}}

This concerns the transient case ($\beta >1$) in the continuous-time setup:
we get 
\begin{equation*}
\overline{\tau }_{x}^{\left( d\right) }:=\left( \overline{\tau }%
_{x,0}+\sum_{k=1}^{\mathcal{N}_{e}}\overline{\tau }_{0,0}^{\left( k\right)
}\right) \mathbf{1}\left( \tau _{x,0}<\infty \right)
\end{equation*}
After $t>\overline{\tau }_{x}^{\left( d\right) }$, the process $\overline{X}%
\left( t\right) $ is a pure-birth Yule process at rate $r_{\overline{X}%
\left( t\right) }$ started at $0.$ A pure birth Yule process started at $%
Y\left( 0\right) $ is defined by

\begin{eqnarray*}
Y\left( t\right) &=&Y\left( 0\right) +P\left( \int_{0}^{t}r_{Y\left(
s\right) }ds\right) \\
Y\left( t\right) &=&Y\left( 0\right) +\int_{0}^{t}r_{Y\left( s\right) }ds+%
\overline{P}\left( \int_{0}^{t}r_{Y\left( s\right) }ds\right)
\end{eqnarray*}
where $P\left( t\right) $ is the standard Poisson process with rate $t$ and $%
\overline{P}\left( t\right) $ the compensated Poisson process (as a
martingale). This process has generator 
\begin{equation*}
G_{Y}h\left( y\right) =r_{y}\left( h\left( y+1\right) -h\left( y\right)
\right) ,
\end{equation*}
meaning 
\begin{equation*}
\Bbb{E}_{y}h\left( Y\left( t\right) \right) =\Bbb{E}_{y}h\left( Y\left(
0\right) \right) +\int_{0}^{t}\Bbb{E}_{y}\left( G_{Y}h\right) \left( Y\left(
s\right) \right) ds.
\end{equation*}
When $h\left( y\right) =y$, with $y\left( t\right) :=\Bbb{E}_{y}Y\left(
t\right) $, $y\left( 0\right) =y$, by Jensen inequality 
\begin{equation*}
\overset{.}{y}\left( t\right) =\Bbb{E}_{y}r\left( Y\left( t\right) \right)
\geq r_{0}\left( y\left( t\right) +1\right) ^{\lambda }\text{ if }\lambda >1,
\end{equation*}
which explodes in finite time.

When the CTMC is transient, it returns to $0$ finitely many times before
drifting to $\infty $ at time $\overline{\tau }_{x}^{\left( d\right) }$ (the
time to its last visit to $0$). From this instant on, it behaves like a pure
birth Yule process: this process grows by one unit and the times separating
these successive unit increments are Exp$(r_{x})$ distributed where $r_{x}$
is the rate at which the transition $x\rightarrow x+1$ proceeds. And here we
chose $r_{x}=r_{0}\left( x+1\right) ^{\lambda }$. Therefore, under the
condition $\lambda >1$ the CTMC explodes in finite time (because then: $%
\sum_{x\geq 0}1/r_{x}<\infty $) with probability $1$. Under this condition,
the final time at which explosion occurs is thus: $\overline{\tau }%
_{x}^{\left( d\right) }+\sum_{x\geq 0}E_{x}.$ If $\lambda \leq 1$, the chain
drifts to $\infty $ without exploding in finite time.

\subsection{\textbf{Continuous state-space}}

Some of the results for the discrete-space, continuous-time, model $%
\overline{X}$ have analogues for a continuous-time Markov chain on the
non-negative reals.

Suppose that given the process $X$ moves up and is in a small neighborhood
of state $x>0$, it moves to a location $y>x$ with probability $f\left(
y-x\right) dx$ where $f$ is a probability density on $\Bbb{R}_{+}$. The
continuous-time, continuous-state version of the process $\overline{X}$, say 
$\underline{X}$, has backward generator $G_{\underline{X}}$ whose action on
real-valued bounded functions $h$ on $\left[ 0,\infty \right) $ now is 
\begin{equation*}
G_{\underline{X}}h\left( x\right) =r_{x}\left\{ \left( p_{x}\int_{0}^{\infty
}\left( h\left( x+z\right) -h\left( x\right) \right) f\left( z\right)
dz+q_{x}\left( h\left( 0\right) -h\left( x\right) \right) \right) \right\} .
\end{equation*}
Whenever it exists, the invariant density measure $\left( \underline{\pi }%
_{x},x\geq 0\right) $ of $\underline{X}$ obeys 
\begin{equation*}
r_{y}\underline{\pi }_{y}=\int_{0}^{y}r_{x}p_{x}\underline{\pi }_{x}f\left(
y-x\right) dx+\delta _{y,0}\int_{0}^{\infty }r_{x}p_{x}\underline{\pi }%
_{x}dx,
\end{equation*}
which is the continuous-state version of the equation $\overline{\mathbf{\pi 
}}^{\prime }Q=\mathbf{0}^{\prime }$ characterizing $\overline{\pi }_{x}$, $%
x\in \Bbb{N}_{0}.$

\section{A variant of the model: the Pareto-Zipf paradigm}

Let $\alpha ,\beta >0.$ Consider the new discrete time-homogeneous Markov
chain (MC) $X:=\left( X_{n};n\geq 0\right) $ with state-space $\Bbb{N}_{0}$
and non-homogeneous spatial transition probabilities now characterized by:

$\bullet $ given $X_{n}=x\in \left\{ 1,2,...\right\} $, the increment of $%
X_{n}$ is 
\begin{equation}
\begin{array}{l}
+1\text{ with probability}:\text{ }p_{x}=\left( 1+x^{-\beta }\right)
^{-\alpha } \\ 
-x\text{ with probability}:\text{ }q_{x}=1-\left( 1+x^{-\beta }\right)
^{-\alpha }.
\end{array}
\label{0c}
\end{equation}

$\bullet $ given $X_{n}=0$, the increment of $X_{n}$ is $+1$ with
probability $p_{0}\leq 1$ and $0$ with probability $q_{0}=1-p_{0}.$

For large $x$, the probabilities $p_{x}$ and $q_{x}$ have a behaviour very
similar to the ones in the previous case (\ref{0a}). We will deal with this
new DTMC, the continuous-time version of which being left to the reader. We
could think of adding a base parameter $\nu >-1$ (giving $p_{x}=\left( \nu
+x^{-\beta }\right) ^{-\alpha }$) but we will not detail this model as it
would deal with the Hurwitz zeta function instead of the zeta function that
we will meet next when $\nu =1$ as in (\ref{0c}). Mutatis mutandis,
proceeding as before, we conclude:

- $C_{1}=\sum_{y=1}^{\infty }q_{y}<\infty $ if and only if $\beta >1$: in
this case the MC is transient with no invariant measure.

- If $\beta <1,$ then $C_{1}=\infty $ and $C_{2}=\sum_{x\geq
1}\prod_{y=0}^{x-1}p_{y}<\infty :$ the MC is positive recurrent.
Furthermore, for large $x$: 
\begin{equation*}
\pi _{x}\sim \prod_{y=0}^{x-1}\left( 1+y^{-\beta }\right) ^{-\alpha }\sim
e^{-\alpha \sum^{x}y^{-\beta }}\sim e^{-\frac{\alpha }{1-\beta }x^{1-\beta }}
\end{equation*}
with stretched exponential behaviour.

- In the critical case $\beta =1$: 
\begin{equation*}
\pi _{x}\sim \prod_{y=0}^{x-1}\left( 1+y^{-1}\right) ^{-\alpha }\sim
e^{-\alpha \sum^{x}y^{-1}}\sim x^{-\alpha }
\end{equation*}
with power-law$\left( \alpha \right) $ behaviour. The chain is recurrent. We
will compute the exact shape of $\pi _{x}$. But from this one concludes that
the DT chain is positive recurrent if $\alpha >1$, null recurrent if $\alpha
\leq 1.$ To make this precise, let us consider the

\subsection{\textbf{Return time to the origin.}}

They are the times $\tau _{0,0}\geq 1$ between consecutive visits to $0$
(the first return times to $0$). We have 
\begin{equation}
\Bbb{P}\left( \tau _{0,0}=x+1\right) =q_{x}\prod_{y=0}^{x-1}p_{y}\text{, }%
x\geq 0.  \label{13d}
\end{equation}
Equivalently, $\Bbb{P}\left( \tau _{0,0}>x\right) =\prod_{y=0}^{x-1}p_{y}.$
Recall that $\Bbb{P}\left( \tau _{0,0}<\infty \right) =1$ if and only if $%
C_{1}=\infty $ which is the recurrence condition for $X_{n}$.

We conclude that: If $C_{1}\equiv \sum_{y=1}^{\infty }q_{y}<\infty ,$ $X_{n}$
is transient with $\Bbb{P}\left( \tau _{0,0}<\infty \right) <1$. If $%
C_{1}=\infty ,$ $X_{n}$ is recurrent. If $C_{1}=\infty $ and $C_{2}\equiv
\sum_{x\geq 1}\prod_{y=0}^{x-1}p_{y}<\infty ,$ $X_{n}$ is positive recurrent
with $\mu :=\Bbb{E}\left( \tau _{0,0}\right) =1/\pi _{0}=1+C_{2}<\infty .$
If $C_{1}=C_{2}=\infty ,$ $X_{n}$ is null recurrent with $\tau _{0,0}<\infty 
$ almost surely (a.s.) and $\Bbb{E}\left( \tau _{0,0}\right) =\infty $.
Recall that in the positive recurrent case, 
\begin{equation*}
\pi _{x}=\frac{\prod_{y=0}^{x-1}p_{y}}{1+C_{2}}=\frac{\Bbb{P}\left( \tau
_{0,0}>x\right) }{\mu }\text{, }x\geq 0.
\end{equation*}
\newline

APPLICATION:

- If $\beta >1:$ the chain is transient with ($0<u_{\infty }<1$) 
\begin{equation*}
\Bbb{P}\left( \tau _{0,0}=\infty \right) =u_{\infty }=\prod_{y=0}^{\infty
}\left( 1+y^{-\beta }\right) ^{-\alpha }>0.
\end{equation*}

- If $\beta <1:$ the chain is positive recurrent with 
\begin{equation*}
\mu :=\Bbb{E}\left( \tau _{0,0}\right) =1+C_{2}=1+\sum_{x\geq
1}\prod_{y=0}^{x-1}\left( 1+y^{-\beta }\right) ^{-\alpha }<\infty .
\end{equation*}

- Critical case: $\beta =1$, $\alpha >0$. In this case, let 
\begin{eqnarray*}
\psi _{0}\left( z\right) &:&=\sum_{x\geq
1}q_{x}z^{x}\prod_{y=1}^{x-1}p_{y}=\sum_{x\geq 1}z^{x}\left( 1-\left(
1+x^{-1}\right) ^{-\alpha }\right) \prod_{y=1}^{x-1}\left( 1+y^{-1}\right)
^{-\alpha } \\
&=&\sum_{x\geq 1}z^{x}\left( 1-\left( 1+x^{-1}\right) ^{-\alpha }\right)
x^{-\alpha }=1-\frac{1}{z}\left( 1-z\right) \text{Li}_{\alpha }\left(
z\right)
\end{eqnarray*}
involving the polylog function: Li$_{\alpha }\left( z\right) :=\sum_{x\geq
1}z^{x}x^{-\alpha }$, obeying 
\begin{equation*}
\text{Li}_{\alpha }\left( z\right) =\int_{0}^{z}\frac{\text{Li}_{\alpha
-1}\left( z^{\prime }\right) }{z^{\prime }}dz^{\prime }.
\end{equation*}

The function $\psi _{0}\left( z\right) $ is the probability generating
function (pgf) of a discrete Pareto rv $P\geq 1$ with tail parameter $\alpha 
$ \cite{BK} ($\Bbb{P}\left( P>x\right) =\left( x+1\right) ^{-\alpha }$, $%
x\in \Bbb{N}_{0}$). In our setup, the rv $P$ can be identified to $\tau
_{0,0}^{+}-1\overset{d}{=}H^{+}$. The polylogarithm can be expressed in
terms of the integral of the Bose-Einstein distribution

\begin{equation*}
\text{Li}_{\alpha }\left( z\right) =\frac{1}{\Gamma \left( \alpha \right) }%
\int_{0}^{\infty }\frac{x^{\alpha -1}}{z^{-1}e^{x}-1}dx=\frac{z}{\Gamma
\left( \alpha \right) }\int_{0}^{1}\frac{\left( -\log u\right) ^{\alpha -1}}{%
1-uz}du,
\end{equation*}
showing that 
\begin{equation*}
\Bbb{P}\left( P>x\right) =\left( x+1\right) ^{-\alpha }=\int_{0}^{1}u^{x}\pi
\left( du\right) \text{ where }\pi \left( du\right) =\frac{1}{\Gamma \left(
\alpha \right) }\left( -\log u\right) ^{\alpha -1}du\text{ }
\end{equation*}
is the probability density of $U=e^{-X}$, with $X\overset{d}{\sim }$Gamma$%
\left( \alpha ,1\right) .$\ The laws of $P\geq 1$\ and $P-1$\ are completely
monotone (so $P-1$\ is ID). We don't know if $P-1$ (with pgf $z^{-1}\psi
_{0}\left( z\right) $) is SD. See the Appendix $8.1$. Note the simplified
expression 
\begin{equation*}
\psi _{0}\left( z\right) =\frac{z}{\Gamma \left( \alpha \right) }\int_{0}^{1}%
\frac{1}{1-uz}\nu \left( du\right) \text{ where }\nu \left( du\right)
=\left( 1-u\right) \pi \left( du\right) .
\end{equation*}

In our setup, 
\begin{equation*}
C_{2}\equiv \sum_{x\geq 1}\prod_{y=0}^{x-1}p_{y}=p_{0}\sum_{x\geq
1}x^{-\alpha }=p_{0}\zeta \left( \alpha \right)
\end{equation*}
involving the zeta function, which is finite only when $\alpha >1$ (positive
recurrence criterion of $X$).

When $\alpha =1$, an easy computation shows that

\begin{equation*}
\psi _{0}\left( z\right) =1+\frac{1}{z}\left( 1-z\right) \log \left(
1-z\right) .
\end{equation*}
We have: 
\begin{eqnarray*}
\Bbb{P}\left( P=x\right) &:&=\rho _{x}=\left[ z^{x}\right] \psi _{0}\left(
z\right) \underset{x\rightarrow \infty }{\sim }\alpha x^{-\left( \alpha
+1\right) }\text{ if }\alpha \neq 1 \\
\Bbb{P}\left( P>x\right) &=&\frac{1}{1+x}\text{ if }\alpha =1\text{ }
\end{eqnarray*}
and for all $x\geq 1$%
\begin{equation}
\rho _{x+1}/\rho _{x}<1,\text{ }\rho _{1}=1-2^{-\alpha }  \label{f13b}
\end{equation}
($\rho _{x}$ is monotone decreasing). As a result, $\rho _{1}$ is the
maximal value of the $\rho _{x}$: the probability mass function (pmf) $\rho
_{x}$ has its mode at $x=1.$

Equivalently,\emph{\ } 
\begin{equation*}
\begin{array}{l}
\psi _{0}\left( z\right) \underset{z\downarrow 1}{\sim }1-\Gamma \left(
1-\alpha \right) \left( 1-z\right) ^{\alpha }\text{ if }0<\alpha <1\text{ }
\\ 
\psi _{0}\left( z\right) \underset{z\downarrow 1}{\sim }1-\zeta \left(
\alpha \right) \left( 1-z\right) +O\left( 1-z\right) ^{\alpha }\text{ if }%
\alpha >1 \\ 
\psi _{0}\left( z\right) \underset{z\downarrow 1}{\sim }1+\left( 1-z\right)
\log \left( 1-z\right) \text{ if }\alpha =1
\end{array}
\end{equation*}

\subsection{\textbf{Return time to }$0$}

The pgf of $\tau _{0,0}$ itself therefore reads,

\begin{eqnarray*}
\phi _{0,0}\left( z\right) &=&\sum_{x\geq
0}q_{x}z^{x+1}\prod_{y=0}^{x-1}p_{y}=z\left( q_{0}+p_{0}\psi _{0}\left(
z\right) \right) \\
&=&z\left( 1-\frac{p_{0}}{z}\left( 1-z\right) \text{Li}_{\alpha }\left(
z\right) \right) .
\end{eqnarray*}
We find the mean persistence time 
\begin{equation*}
\phi _{0,0}^{\prime }\left( 1\right) =\Bbb{E}\left( \tau _{0,0}\right) :=\mu
=1+C_{2}=\frac{1}{\pi _{0}}=1+p_{0}\zeta \left( \alpha \right) .
\end{equation*}
It can be checked that if $\Bbb{E}\left( \tau _{0,0}\right) $\ is to exist,
then necessarily $\Bbb{E}\left( \tau _{0,0}\right) >1+p_{0}$.

* This condition forces $\alpha >1$\ which is the positive recurrence
condition for the critical MCC. If $p_{0}=1$, $\Bbb{E}\left( \tau
_{0,0}\right) =1+\zeta \left( \alpha \right) >2$ and there are no trivial
excursions.

\subsection{\textbf{Contact probability at }$0$.}

With $\phi _{0,0}\left( z\right) =z\left( q_{0}+p_{0}\psi _{0}\left(
z\right) \right) $ in the critical special case.

- If $0<\alpha <1$, using $\psi _{0}\left( z\right) \underset{z\downarrow 1}{%
\sim }1-\Gamma \left( 1-\alpha \right) \left( 1-z\right) ^{\alpha }$, we
have 
\begin{equation*}
1-\phi _{0,0}\left( z\right) \underset{z\downarrow 1}{\sim }p_{0}\left(
1-\psi _{0}\left( z\right) \right) \sim p_{0}\Gamma \left( 1-\alpha \right)
\left( 1-z\right) ^{\alpha }
\end{equation*}
showing by singularity analysis that, in the null recurrent case (algebraic
decay of the contact probability) 
\begin{equation*}
\Bbb{P}_{0}\left( X_{n}=0\right) \underset{n\uparrow \infty }{\sim }\frac{1}{%
p_{0}\Gamma \left( 1-\alpha \right) }n^{-\left( 1-\alpha \right) }.
\end{equation*}
When $\alpha =1-\varepsilon $ ($\varepsilon >0$ small), the constant in
front of $n^{-\left( 1-\alpha \right) }$ vanishes like $\varepsilon /p_{0}$.

- When $\alpha =1$, using $\psi _{0}\left( z\right) \underset{z\downarrow 1}{%
\sim }1+\left( 1-z\right) \log \left( 1-z\right) $%
\begin{equation*}
1-\phi _{0,0}\left( z\right) =1-z\left( q_{0}+p_{0}\psi _{0}\left( z\right)
\right) \underset{z\downarrow 1}{\sim }p_{0}\left( 1-z\right) \log \left(
1-z\right) ,
\end{equation*}
with logarithmic singularity, showing by singularity analysis, that 
\begin{equation*}
\Bbb{P}_{0}\left( X_{n}=0\right) \underset{n\uparrow \infty }{\sim }\frac{1}{%
p_{0}\log n}.
\end{equation*}

- When $\alpha >1$ (positive recurrence case), using $\psi _{0}\left(
z\right) \underset{z\downarrow 1}{\sim }1-\zeta \left( \alpha \right) \left(
1-z\right) +O\left( 1-z\right) ^{\alpha },$ 
\begin{equation*}
1-\phi _{0,0}\left( z\right) \underset{z\downarrow 1}{\sim }\left(
1-z\right) \left( 1+p_{0}\zeta \left( \alpha \right) \right)
\end{equation*}
showing, as required that $\Bbb{P}_{0}\left( X_{n}=0\right) \underset{%
n\uparrow \infty }{\rightarrow }\pi _{0}=1/\left( 1+p_{0}\zeta \left( \alpha
\right) \right) .$ When $\alpha =1+\varepsilon $ ($\varepsilon >0$ small),
using $\zeta \left( \alpha \right) \underset{\alpha \uparrow 1}{\sim }%
\varepsilon ^{-1}$: $\pi _{0}\sim \varepsilon /p_{0}$, just like when $%
\alpha <1.$

\subsection{\textbf{Invariant measure }$\left( \alpha >1\right) $\textbf{.}}

In the positive recurrent case, we also have

\begin{eqnarray*}
\pi _{x} &=&\frac{\prod_{y=0}^{x-1}p_{y}}{\mu }=\pi
_{0}p_{0}\prod_{y=1}^{x-1}\left( 1+y^{-1}\right) ^{-\alpha }=\pi
_{0}p_{0}x^{-\alpha }\text{, }x\geq 0 \\
\frac{\pi _{x}}{1-\pi _{0}} &=&\Bbb{P}\left( X_{\infty }=x\mid X_{\infty
}\geq 1\right) =\frac{x^{-\alpha }}{\zeta \left( \alpha \right) }\text{, }%
x\geq 1
\end{eqnarray*}
an identity showing that $Y_{\infty }:=X_{\infty }\mid X_{\infty }\geq 1%
\overset{d}{\sim }$Zipf$\left( \alpha -1\right) $\ so with (recalling Li$%
_{\alpha }\left( 1\right) =\zeta \left( \alpha \right) $) 
\begin{eqnarray*}
\Bbb{E}\left( z^{Y_{\infty }}\right) &=&\text{Li}_{\alpha }\left( z\right) /%
\text{Li}_{\alpha }\left( 1\right) :=\psi _{\infty }\left( z\right) \\
\Bbb{E}\left( z^{X_{\infty }}\right) &=&\pi _{0}+\left( 1-\pi _{0}\right)
\psi _{\infty }\left( z\right) =\left( 1+p_{0}\text{Li}_{\alpha }\left(
z\right) \right) /\left( 1+p_{0}\text{Li}_{\alpha }\left( 1\right) \right) .
\end{eqnarray*}
This gives an explicit expression of the pgf of $X_{\infty }$ (in the
positive recurrent case) in terms of a polylog function$.$ As in the
previous example, $X_{\infty }$ is the mixture of the two rvs $Y_{\infty
}\geq 1$ and (say) $Y_{0}$ which is degenerate at $0$: $Y_{0}\overset{d}{%
\sim }\delta _{0}.$

The shifted Zipf rv $Y_{\infty }-1:=\left( X_{\infty }-1\mid X_{\infty }\geq
1\right) \geq 0,$ \cite{New}, with pmf $\zeta \left( \alpha \right)
^{-1}\left( x+1\right) ^{-\alpha }$, $x\geq 0$ and pgf $z^{-1}\psi _{\infty
}\left( z\right) $, is discrete self-decomposable (\cite{SH}, Example $12.18$%
, page $435$)\footnote{%
Caution: What is called here a Zipf rv is called a Pareto rv in \cite{SH}.}.
It has tail index $\alpha -1>0$, with finite mean if $\alpha >2.$

As in the previous model, we can ask under what conditions on $p_{0}$ the rv 
$X_{\infty }$\emph{\ }itself, can be ID or not. Firstly, we observe that $%
\pi _{x}$ is unimodal with mode at the origin ($\pi _{x+1}/\pi _{x}=\left(
1+x^{-1}\right) ^{-\alpha }<1$, $x\geq 1$, $\pi _{1}/\pi _{0}=p_{0}<1$).%
\emph{\ }Observing next that $\pi _{x}^{2}\leq \pi _{x-1}\pi _{x+1}$ for $%
x\geq 2$,\emph{\ }a sufficient condition for\emph{\ }$X_{\infty }$ to be
log-convex and so infinitely divisible is that 
\begin{equation*}
\Bbb{P}\left( X_{\infty }=1\right) ^{2}\leq \Bbb{P}\left( X_{\infty
}=0\right) \Bbb{P}\left( X_{\infty }=2\right) \text{ which is }p_{0}\leq
2^{-\alpha }.
\end{equation*}
If $f$ is an unimodal continuous probability density, the discrete pmf 
\begin{equation*}
\pi _{x}=\frac{1}{x!}\int_{0}^{\infty }y^{x}e^{-y}f\left( y\right) dy\text{, 
}x\in \Bbb{N}_{0}
\end{equation*}
obtained by a Poisson--intensity-mixing procedure is infinitely divisible
and unimodal, \cite{Hol}. We don't know if there are values of $p_{0},$
small enough (else $\pi _{0}$ large enough), for which $f$ is SD and
unimodal, so that\emph{\ }$X_{\infty }$\emph{\ }would be
discrete-self-decomposable (see the Appendix).

Note finally that, 
\begin{equation*}
\Bbb{E}\left( Y_{\infty }^{q}\right) =\frac{1}{\zeta \left( \alpha \right) }%
\sum_{x\geq 1}\frac{1}{x^{\alpha -q}}=\frac{\zeta \left( \alpha -q\right) }{%
\zeta \left( \alpha \right) },\text{ }0\leq q<\alpha -1.
\end{equation*}
Owing to the Euler product formula for zeta function over the primes $%
\mathcal{P}=\left\{ 2,3,5,...\right\} $%
\begin{equation*}
\zeta \left( \alpha \right) =\prod_{p\in \mathcal{P}}\left( 1-p^{-\alpha
}\right) ^{-1},
\end{equation*}
this shows that 
\begin{equation*}
Y_{\infty }\overset{d}{=}\prod_{p\in \mathcal{P}}p^{G\left( p^{-\alpha
}\right) }\text{,}
\end{equation*}
where $G\left( p^{-\alpha }\right) $, $p\in \mathcal{P}$ is a sequence of
independent Geometric rvs with $\Bbb{P}\left( G\left( p^{-\alpha }\right)
=x\right) =p^{-\alpha x}\left( 1-p^{-\alpha }\right) $, $x\geq 0$. We note
therefore that the pgf of 
\begin{equation*}
\log \left( Y_{\infty }\right) =\sum_{p\in \mathcal{P}}G\left( p^{-\alpha
}\right) \log p
\end{equation*}
is 
\begin{equation*}
\prod_{p\in \mathcal{P}}\frac{1-p^{-\alpha }}{1-z^{\log p}p^{-\alpha }},
\end{equation*}
an infinite product of Geometric$\left( p^{-\alpha }\right) $ pgfs on
multiples of $\log p$, \cite{Lin}. Each such rv is self-decomposable and so
is the infinite product: the rv $\log \left( Y_{\infty }\right) $ is
self-decomposable, an additional remarkable divisibility feature of $%
Y_{\infty }$. The moment function of $X_{\infty }$ itself is 
\begin{eqnarray*}
\Bbb{E}\left( X_{\infty }^{q}\right) &=&\pi _{0}\delta _{q,0}+\left( 1-\pi
_{0}\right) \frac{\zeta \left( \alpha -q\right) }{\zeta \left( \alpha
\right) } \\
&=&\pi _{0}\delta _{q,0}+\frac{p_{0}\zeta \left( \alpha -q\right) }{%
1+p_{0}\zeta \left( \alpha \right) }\text{, }0\leq q<\alpha -1.
\end{eqnarray*}
\newline

* If $\alpha \leq 1$, the special MCC is null recurrent with $\Bbb{E}\left(
\tau _{0,0}\right) =\infty $.

\subsection{\textbf{Renewal aspects.}}

In the positive recurrent case, we can also reproduce the arguments related
to renewal theory, for which 
\begin{equation*}
\Delta =\tau _{0,0}^{0}+\tau _{0,0}^{+},
\end{equation*}
where $\tau _{0,0}^{0}\overset{d}{\sim }$Geo$\left( p_{0}\right) $\ (viz $%
\Bbb{P}\left( \tau _{0,0}^{0}=k\right) =p_{0}q_{0}^{k}$, $k\geq 0$) is the
idle period (the number of consecutive trivial excursions) and $\tau
_{0,0}^{+}\geq 2$\ the length of a true excursion with: $\Bbb{P}\left( \tau
_{0,0}^{+}=x+1\right) =q_{x}\prod_{y=1}^{x-1}p_{y}$, $x\geq 1$, the lengths $%
\tau _{0,0}^{0}$, $\tau _{0,0}^{+}$ being mutually independent.

We observe that the random variable $\tau _{0,0}^{+}-1$ (also the height $%
H^{+}$ of a true excursion), with support $\Bbb{N,}$ can be generated
according to the success run procedure: 
\begin{equation}
\tau _{0,0}^{+}-1\overset{d}{=}H^{+}=\inf \left( x\geq 1:B_{x}\left( \alpha
\right) =1\right) ,  \label{heig}
\end{equation}
where $\left( B_{x}\left( \alpha \right) ;\text{ }x\geq 1\right) $ is a
sequence of independent Bernoulli rvs obeying $\Bbb{P}\left( B_{x}\left(
\alpha \right) =1\right) =1-\left( 1+x^{-1}\right) ^{-\alpha }=q_{x}.$ It is
the Pareto rv $P$ with pgf $\psi _{0}\left( z\right) $ defined above.\newline

APPLICATION (positive recurrent case, $\beta =1$, $\alpha >1$): Clearly,
with $C_{2}^{+}=\sum_{x\geq 1}\prod_{y=1}^{x-1}p_{y}=C_{2}/p_{0}=\zeta
\left( \alpha \right) <\infty $%
\begin{equation*}
\Bbb{E}\left( \tau _{0,0}^{0}\right) =p_{0}/q_{0}\text{ and }\Bbb{E}\left(
\tau _{0,0}^{+}\right) =:\mu ^{+}=1+C_{2}^{+}=1+\zeta \left( \alpha \right)
>2.
\end{equation*}
gives the average contribution to $\Delta $ of the two components in the
special critical MC case. The $\Delta $s now constitute the interval lengths
of an alternating renewal process.

\section{Appendix}

\subsection{Discrete infinite-divisibility and self-decomposability}

For the sake of completeness, let us briefly introduce the notion of
self-decomposability.

\begin{definition}
Let $X\geq 0$ be an integer-valued random variable. The probability
generating function (pgf) $\phi \left( z\right) :=\Bbb{E}\left( z^{X}\right) 
$ is the one of a discrete self-decomposable (SD) variable $X$ if for any $%
u\in \left( 0,1\right) $, there is a pgf $\phi _{u}\left( z\right) $
(depending on $u$) such that (see \cite{Steu}), 
\begin{equation}
\phi \left( z\right) =\phi \left( 1-u\left( 1-z\right) \right) \cdot \phi
_{u}\left( z\right) .  \label{f32}
\end{equation}
\end{definition}

Define the $u$-thinned version of $X,$ say $u\circ X,$\ as the random sum 
\begin{equation}
u\circ X\overset{d}{=}\sum_{x=1}^{X}B_{x}\left( u\right) ,  \label{g2}
\end{equation}
with $\left( B_{x}\left( u\right) \right) _{x\geq 1}$\ a sequence of iid
Bernoulli variables such that $\Bbb{P}\left( B_{x}\left( u\right) =1\right)
=u$, independent of $X.$

This binomial thinning operator, acting on discrete rvs, has been defined in 
\cite{Steu}; it stands as the discrete version of the change of scale $%
X\rightarrow u\cdot X$ for continuous rvs $X$. If $\phi \left( z\right) $\
is the pgf of the SD random variable $X$ obeying Eq. (\ref{f32}), then $X$\
can be additively (self-)decomposed as 
\begin{equation}
X\overset{d}{=}u\circ X^{\prime }+X_{u}.  \label{f34}
\end{equation}
Here, $X$\ and $X^{\prime }$ have the same distribution and $u\circ
X^{\prime }$\ is independent of the remaining random variable $X_{u}$ with
pgf, say $\phi _{u}\left( z\right) .$ \newline

\textbf{Remark:} The self-decomposability idea also (pre-)exists for
continuous rvs on $\Bbb{R}_{+}$: $X>0$ is said to be SD if, for any $u\in
\left( 0,1\right) $ 
\begin{equation}
X\overset{d}{=}u\cdot X^{\prime }+X_{u}.  \label{f34b}
\end{equation}
with $X$\ and $X^{\prime }$ having the same distribution and $u\cdot
X^{\prime }$\ being independent of the remaining random variable $X_{u}>0.$
If $\Phi \left( \lambda \right) =\Bbb{E}\left( e^{-\lambda X}\right) $, $%
\lambda \geq 0$, is to be the Laplace-Stieltjes transform (LST) of $X$, SD,
then for any $u\in \left( 0,1\right) $, there is a LST $\Phi _{u}\left(
\lambda \right) $ (depending on $u$) such that

\begin{equation*}
\Phi \left( \lambda \right) =\Phi \left( \lambda u\right) \cdot \Phi
_{u}\left( \lambda \right) .
\end{equation*}
The two notions of self-decomposability are related as follows: let $Y>0$ be
a continuous rv. Then $Y$ is self-decomposable if and only if the discrete
random variable supported by $\Bbb{N}_{0}$ defined by: $X=P\left( Y\right) $
(where $P\left( Y\right) $ is a Poisson rv with random intensity $Y$) is
discrete self-decomposable (Corollary $1$ of \cite{Sap}). Indeed, 
\begin{equation*}
\phi _{X}\left( z\right) =\Bbb{E}\left( z^{X}\right) =\Phi _{Y}\left(
1-z\right)
\end{equation*}
and, with $\phi _{X_{u}}\left( z\right) =\Bbb{E}\left( z^{X_{u}}\right)
=\Phi _{u}\left( 1-z\right) =\Bbb{E}\left( e^{-\left( 1-z\right)
Y_{u}}\right) $, the pgf of $X_{u}=P\left( Y_{u}\right) $%
\begin{equation*}
\phi _{X}\left( z\right) =\phi _{X}\left( 1-u\left( 1-z\right) \right) \cdot
\phi _{X_{u}}\left( z\right) \Leftrightarrow \Phi _{Y}\left( \lambda \right)
=\Phi _{Y}\left( \lambda u\right) \cdot \Phi _{Y_{u}}\left( \lambda \right) .%
\text{ }\Box
\end{equation*}
In such cases, the pmf of $X$ is related to the density $f$ of $Y$ by: 
\begin{equation*}
\Bbb{P}\left( X=x\right) =\frac{1}{x!}\int_{0}^{\infty }y^{x}e^{-y}f\left(
y\right) dy,\text{ }x\in \Bbb{N}_{0}.
\end{equation*}
This raises the question of which $Y$ are SD? We first need to recall the
notion of an HCM rv. Following (\cite{SH}, $(5.15)$ p. $371$), a HCM
positive rv is one whose density $f$ obeys that 
\begin{equation*}
\forall y>0:\text{ the function }x\rightarrow f\left( xy\right) f\left(
y/x\right) \text{, }x>0
\end{equation*}
is completely monotone on $\left( 2,\infty \right) $ as a function of $%
z=x+1/x$. Completely monotone functions $h$ obey: $\left( -1\right)
^{k}h^{\left( k\right) }\left( x\right) \geq 0$ for all $k\geq 0$ in some
range of $x$.

With $a>0$, consider now the rv $Y=G\left( a\right) ^{1/\beta }$, the $%
1/\beta -$power of $G\left( a\right) \overset{d}{\sim }$Gamma$\left(
a,1\right) $. This rv is hyperbolically completely monotone (HCM) if and
only if $\left| \beta \right| \leq 1$ (ex. $12.8$ of \cite{SH}). This is
also true of the so-called Generalized Inverse Gaussian rvs \cite{Hal}, with
density 
\begin{equation}
f_{\lambda ,\delta _{1},\delta _{2}}\left( x\right) =\left( \frac{\delta _{1}%
}{\delta _{2}}\right) ^{\alpha /2}\frac{1}{2K_{\alpha }\left( \sqrt{\delta
_{1}\delta _{2}}\right) }x^{\alpha -1}\exp \left\{ -\frac{1}{2}\left( \delta
_{1}x+\delta _{2}/x\right) \right\} \text{, }x>0  \label{eq51c}
\end{equation}
in the parameter range: $\delta _{2}\geq 0$, $\delta _{1}>0$ if $\lambda >0$%
, $\delta _{2}>0$, $\delta _{1}\geq 0$ if $\lambda <0$ and $\delta _{2}>0$, $%
\delta _{1}>0$ if $\lambda =0.$ Such densities include Gamma distributions ($%
\delta _{2}=0$) and inverse Gamma distributions ($\delta _{1}=0$) as
particular cases. $K_{\alpha }$ is the modified Bessel function of the
second kind.

HCM rvs constitute a subclass of Generalized-Gamma-Convolution (GGC) rvs
(Proposition $2$ of \cite{Bondesson}) and GGC rvs are SD (Theorem $1$ of 
\cite{Bondesson}). We refer to \cite{SH}, Section $5$ and \cite{Bond}, for
the precise definition of GGC rvs. So, whenever $Y$ is a GGC rv, it is SD
and $X=P\left( Y\right) $, as a Poisson-mixture with respect to a SD
distribution, is discrete-SD. Self-decomposable distributions are unimodal.%
\newline

Coming back to discrete self-decomposability itself, the following
representation result is also known to hold true, see \cite{Steu}. Let $%
R\left( z\right) $\ (with $r_{0}=R\left( 0\right) >0$) be the canonical
function defined through 
\begin{equation}
\phi \left( z\right) =\Bbb{E}\left( z^{X}\right) =e^{-\int_{z}^{1}R\left(
z^{\prime }\right) dz^{\prime }}.  \label{g3}
\end{equation}
The random variable $X$\ is discrete SD if and only if the function $h\left(
z\right) :=1-\left( 1-z\right) R\left( z\right) /r_{0}$\ defines a pgf such
that $h\left( 0\right) =0$ (see \cite{Sch}, Lemma $2.13$). Consequently, $X$%
\ is discrete SD if and only if, for some $r_{0}>0$, its pgf can be written
in the form 
\begin{equation}
\phi \left( z\right) =e^{-r_{0}\int_{z}^{1}\frac{1-h\left( z^{\prime
}\right) }{1-z^{\prime }}dz^{\prime }}.  \label{g4}
\end{equation}
This means that the series coefficients $\left( r_{x}=\left[ z^{x}\right]
R\left( z\right) ,\text{ }x\in \Bbb{N}_{0}\right) $ of $R\left( z\right) $
constitute a non-negative, non-increasing sequence of $x$ (Theorem $4.13$ p. 
$271$ of \cite{Steu}). As a result, the associated probability system $\Bbb{P%
}\left( X=x\right) :=\pi _{x}$, $x\in \Bbb{N}_{0}$ of $X$, if SD, is
unimodal, with mode at the origin if and only if $r_{0}=\frac{\pi _{1}}{\pi
_{0}}\leq 1$. The SD subclass of infinitely divisible distributions (ID)
therefore consists of unimodal distributions, with mode possibly at the
origin (Theorem $2.3$ of \cite{Steu}).

Note that $X$ is ID if and only if the sequence $\left( r_{x},x\in \Bbb{N}%
_{0}\right) $ is non-negative only, with (as can be checked from (\ref{g3}))
the sequences $\left( \pi _{x},\text{ }r_{x},\text{ }x\in \Bbb{N}_{0}\right) 
$ related by the convolution formula 
\begin{equation*}
\left( x+1\right) \pi _{x}=\sum_{y=0}^{x}\pi _{y}r_{x-y}\text{, }x\in \Bbb{N}%
_{0}.
\end{equation*}
We recall that for rvs with integral support $\Bbb{N}_{0},$ the notion of an
infinitely divisible rv coincides with the one of a compound Poisson rv.

\begin{definition}
A compound Poisson rv is one which is obtained as an independent Poisson sum
of positive iid rvs which are called compounding rvs.
\end{definition}

A rv with probability mass $0$ at $0$ cannot be ID.

With $r>0$, the pgf of ID rvs takes the form: $\phi \left( z\right)
=e^{-r\left( 1-h\left( z\right) \right) }$ where $h\left( z\right) $ is the
pgf of the compounding rvs, obeying $h\left( 0\right) =0$.\newline

\textbf{Remark}: Let $Y_{\infty }\geq 1$ be a rv such that $Y_{\infty }-1$
is SD. In the study of the MCC we encountered the delicate problem of
deciding whether or not the mixed rv $X_{\infty }\overset{d}{\sim }\pi
_{0}\delta _{0}+\left( 1-\pi _{0}\right) Y_{\infty }$ was ID or SD. Consider
the simpler case where $Y_{\infty }-1$ is Geometric$\left( p\right) ,$ which
is SD. Then 
\begin{eqnarray*}
\phi _{X_{\infty }}\left( z\right) &=&\pi _{0}+\left( 1-\pi _{0}\right) 
\frac{pz}{1-qz}=\frac{\pi _{0}-\left( \pi _{0}-p\right) z}{1-qz} \\
\log \phi _{X_{\infty }}\left( z\right) &=&\log \left( \pi _{0}-\left( \pi
_{0}-p\right) z\right) -\log \left( 1-qz\right) \\
R\left( z\right) &=&\log \phi _{X_{\infty }}\left( z\right) ^{\prime }=-%
\frac{\pi _{0}-p}{\pi _{0}-\left( \pi _{0}-p\right) z}+\frac{q}{1-qz} \\
r_{x} &=&\left[ z^{x}\right] R\left( z\right) =q^{x+1}-\left( 1-\frac{p}{\pi
_{0}}\right) ^{x+1}
\end{eqnarray*}
Observing $q\geq 1-\frac{p}{\pi _{0}}$, a necessary and sufficient condition
for $X_{\infty }$ to be ID is $1-\frac{p}{\pi _{0}}\geq 0$ ($\pi _{0}\geq p$%
), leading to $r_{x}\geq 0$ for all $x\geq 0$. Under this condition, it will
be SD if and only if $r_{x+1}\leq r_{x}$, meaning 
\begin{equation*}
\left( \frac{1-\frac{p}{\pi _{0}}}{q}\right) ^{x}\leq \pi _{0}\text{ for all 
}x.
\end{equation*}
This will always be the case if in addition $1-\frac{p}{\pi _{0}}\leq q\pi
_{0}$, equivalently if $\pi _{0}\leq p/q.$ If $p\geq q$ ($p\geq 1/2$), $%
X_{\infty }$ is SD if $\pi _{0}\geq p.$ If $p<q$ ($p<1/2$), $X_{\infty }$
will be SD only if $p\leq \pi _{0}\leq p/q$. In the range $\pi _{0}\in
\left( p/q,1\right) $ it is only ID. $\Box $\newline

\textbf{Complete monotonicity:} A rv $X$ with support $\Bbb{N}$ is
completely monotone if, for some probability measure $\pi $ on $\left[
0,1\right] $, the following Hausdorff representation holds 
\begin{eqnarray*}
\overline{F}\left( x\right) &:&=\Bbb{P}\left( X>x\right)
=\int_{0}^{1}u^{x}\pi \left( du\right) \text{, }x\in \left\{
0,1,2,...\right\} \\
\Bbb{P}\left( X=x\right) &=&\pi _{x}=\int_{0}^{1}u^{x}\left( 1-u\right) \pi
\left( du\right) \text{, }x\in \left\{ 1,2,...\right\} .
\end{eqnarray*}
If this is the case, for all $x\in \left\{ 0,1,2,...\right\} $ (see p. $77$
of \cite{SH}) 
\begin{equation*}
\left( -1\right) ^{k}\Delta ^{\left( k\right) }\overline{F}\left( x\right)
\geq 0\text{, equivalently }\left( -1\right) ^{k}\Delta ^{\left( k\right)
}\pi _{x}\geq 0
\end{equation*}
where $\Delta :$ $\Delta h\left( x\right) =h\left( x+1\right) -h\left(
x\right) $ is the right-shift operator and $\Delta ^{\left( k\right) }$ its $%
k-$th iterate. If this is the case, the rv $X-1,$ with support $\Bbb{N}_{0}$%
, is completely monotone, log-convex and therefore infinitely divisible (see
Theorem $10.4$ p. $77$ of \cite{SH}).

\subsection{On simple Markov realizations of ID and SD distributions}

The catastrophe Markov chains that have been studied in this draft showed up
invariant measures that can be either ID or SD (such as $Y_{\infty }-1$
always or $X_{\infty }$ itself in some parameter range on $p_{0}$).
Different Markov processes can have the same invariant equilibrium measure
and here are natural ones whenever the latter is either ID or SD:\newline

- The ID case: Consider a time-inhomogeneous Poisson process $P\left(
R_{t}\right) $ with decaying rate function $re^{-t}$ and intensity $%
R_{t}=r\left( 1-e^{-t}\right) $, $r>0.$ Consider the compound-Poisson
process (with independent increments): 
\begin{equation*}
X_{t}=\sum_{k=1}^{P\left( R_{t}\right) }\Delta _{k}
\end{equation*}
where $\left( \Delta _{k};k\geq 1\right) $ is the iid sequence of the
positive jumps occurring at the jump times of $P\left( R_{t}\right) $. If $%
h\left( z\right) =\Bbb{E}\left( z^{\Delta _{1}}\right) $, 
\begin{equation*}
\phi _{t}\left( z\right) =\Bbb{E}\left( z^{X_{t}}\right) =e^{-R_{t}\left(
1-h\left( z\right) \right) }\underset{t\rightarrow \infty }{\rightarrow }%
e^{-r\left( 1-h\left( z\right) \right) }
\end{equation*}
which is the pgf of an ID rv. A mechanism responsible of the decay of the
population, when balanced by incoming immigrants with sizes $\Delta ,$ will
produce an ID limiting population size.\newline

- The SD case (subcritical branching with immigration): Consider now a
time-homogeneous compound Poisson process $P\left( rt\right) ,$ $t\geq 0$, $%
P\left( 0\right) =0,$ so with pgf 
\begin{equation}
\Bbb{E}\left( z^{P\left( rt\right) }\right) =\exp \left\{ -rt\left(
1-h\left( z\right) \right) \right\} ,  \label{f41}
\end{equation}
where $h\left( z\right) $, with $h\left( 0\right) =0$, is the pgf of the
jumps arriving at the jump times of $P\left( rt\right) $ having rate $r>0.$
Let now 
\begin{equation}
\varphi _{t}\left( z\right) =1-e^{-t}\left( 1-z\right) \text{,}  \label{f41b}
\end{equation}
be the pgf of a pure-death (rate-$1$) Greenwood branching process started
with one particle at $t=0$, \cite{Green}. This expression of $\varphi
_{t}\left( z\right) $ is easily seen to be the solution to $\overset{.}{%
\varphi }_{t}\left( z\right) =f\left( \varphi _{t}\left( z\right) \right)
=1-\varphi _{t}\left( z\right) $, $\varphi _{0}\left( z\right) =z$, as is
usual for a pure-death continuous-time Bellman-Harris branching processes
with affine branching mechanism $f\left( z\right) =r_{d}\left( 1-z\right) $
and fixing the death rate to be $r_{d}=1$, \cite{Harris}. The distribution
function of the lifetime of the initial particle is thus $1-e^{-t}$. Let $%
X_{t}$ with $X_{0}=0$ be a random process counting the current size of some
population for which a random number of individuals (determined by $h\left(
z\right) $) immigrate at the jump times of $P\left( rt\right) ,$ each of
which being independently and immediately subject to the latter pure death
Greenwood process. We have 
\begin{equation}
\phi _{t}\left( z\right) :=\Bbb{E}\left( z^{X_{t}}\right) =\exp
-r\int_{0}^{t}\left( 1-h\left( \varphi _{t-s}\left( z\right) \right) \right)
ds\text{, }\phi _{0}\left( z\right) =1,  \label{f42}
\end{equation}
with $\phi _{t}\left( 0\right) =\Bbb{P}\left( X_{t}=0\right) =\exp
-r\int_{0}^{t}\left( 1-h\left( 1-e^{-s}\right) \right) ds,$ the probability
that the population is extinct at $t$. It holds 
\begin{equation}
\phi _{t}\left( z\right) =e^{-r\int_{0}^{t}\left( 1-h\left( 1-e^{-s}\left(
1-z\right) \right) \right) ds}=e^{-r\int_{\varphi _{0}\left( z\right)
=z}^{\varphi _{t}\left( z\right) }\frac{1-h\left( z^{\prime }\right) }{%
1-z^{\prime }}dz^{\prime }}\underset{t\rightarrow \infty }{\rightarrow }\phi
_{\infty }\left( z\right) =e^{-r\int_{z}^{1}\frac{1-h\left( z^{\prime
}\right) }{1-z^{\prime }}dz^{\prime }}.  \label{f43}
\end{equation}
So, $X:=X_{\infty },$ as the limiting population size of this pure-death
process with immigration, is a SD rv, \cite{vanH}. As in the ID case, a
mechanism responsible of the decay of the population is balanced by incoming
immigrants. \newline

\textbf{Remark:} If instead of a pure-death Greenwood branching process, the
immigrants shrink, more generally, according to any \emph{subcritical}
branching process with branching mechanism $f$, as those whose pgf obeys $%
\overset{.}{\varphi }_{t}\left( z\right) =f\left( \varphi _{t}\left(
z\right) \right) $, then, with $f^{\prime }\left( 1\right) <0,$ $\varphi
_{t}\left( z\right) \rightarrow 1$ as $t\rightarrow \infty $ and 
\begin{equation*}
\phi _{t}\left( z\right) =e^{-r\int_{0}^{t}\left( 1-h\left( \varphi
_{s}\left( z\right) \right) \right) ds}=e^{-r\int_{\varphi _{0}\left(
z\right) =z}^{\varphi _{t}\left( z\right) }\frac{1-h\left( z^{\prime
}\right) }{f\left( z^{\prime }\right) }dz^{\prime }}\rightarrow \phi
_{\infty }\left( z\right) =e^{-r\int_{z}^{1}\frac{1-h\left( z^{\prime
}\right) }{f\left( z^{\prime }\right) }dz^{\prime }}.
\end{equation*}
The obtained limiting pgf is the one of a self-decomposable rv induced by
the subcritical semigroup $\varphi _{t}\left( z\right) $ generated by $%
f\left( z\right) $, \cite{vanH}. Recall $f\left( z\right) =\varphi \left(
z\right) -z$ where $\varphi \left( z\right) $ is the pgf of the branching
number per capita in a Bellman-Harris process. $\Box $\newline

\textbf{Acknowledgments:} T.H. acknowledges partial support from the labex
MME-DII (Mod\`{e}les Math\'{e}matiques et \'{E}conomiques de la Dynamique,
de l' Incertitude et des Interactions), ANR11-LBX-0023-01. This work also
benefited from the support of the Chair ``Mod\'{e}lisation Math\'{e}matique
et Biodiversit\'{e}'' of Veolia-Ecole Polytechnique-MNHN-Fondation X.

\end{document}